\documentclass[10pt]{amsart}

\usepackage{amssymb}
\usepackage[english]{babel}
\usepackage[DIV10, headinclude, BCOR 10mm]{typearea}
\usepackage{enumerate}
\usepackage[colorlinks,linkcolor=blue,citecolor=blue,urlcolor=blue,a4paper]{hyperref}
\usepackage{amsrefs}
\usepackage[all]{xy}

\hypersetup{pdfauthor={Roman Sauer, Andreas Thom},
                pdftitle ={A spectral sequence to compute $L^2$-Betti numbers of groups and groupoids},
           }

\newcommand{\cala}{{\mathcal A}}
\newcommand{\calb}{{\mathcal B}}
\newcommand{\calc}{{\mathcal C}}
\newcommand{\calm}{{\mathcal M}}
\newcommand{\calq}{{\mathcal Q}}
\newcommand{\calg}{{\mathcal G}}
\newcommand{\calr}{{\mathcal R}}
\newcommand{\cals}{{\mathcal S}}
\newcommand{\calt}{{\mathcal T}}
\newcommand{\calu}{{\mathcal U}}

\newcommand{\bbC}{{\mathbb C}}
\newcommand{\bbN}{{\mathbb N}}
\newcommand{\bbR}{{\mathbb R}} 
\newcommand{\bbZ}{{\mathbb Z}} 
\newcommand{\bbQ}{{\mathbb Q}} 

\newcommand{\sauer}[1]%
{\marginpar{\begin{minipage}[c]{\marginparwidth}%
                                  \scriptsize\emph{\color{blue}Roman: #1}%
                                  \end{minipage}}%
}
\newcommand{\andreas}[1]%
{\marginpar{\begin{minipage}[c]{\marginparwidth}%
                                  \scriptsize\emph{\color{red}Andreas: #1}%
                                  \end{minipage}}%
}

\newcommand{\aut}{\operatorname{Aut}}
\newcommand{\emor}{\operatorname{End}}
\newcommand{\betti}{b^{(2)}}
\newcommand{\bs}{\backslash}
\newcommand{\defq}{\mathrel{\mathop:}=}
\newcommand{\colim}{\operatorname{colim}}
\newcommand{\id}{\operatorname{id}}
\renewcommand{\Im}{\operatorname{im}}
\newcommand{\modules}[1]{\operatorname{Mod}(#1)}
\newcommand{\locmod}[1]{\modules{#1}_{\operatorname{loc}}}
\newcommand{\compmod}[1]{\modules{#1}_{\operatorname{comp}}}
\newcommand{\pr}{\operatorname{pr}}

\newcommand{\cd}{\operatorname{cd}}

\newcommand{\Tor}{{\rm Tor}}
\newcommand{\Ext}{{\rm Ext}}
\newcommand{\hyp}{{\rm hyp}}
\newcommand{\SL}{\operatorname{SL}}
\newcommand{\SO}{\operatorname{SO}}
\newcommand{\map}{\operatorname{map}}

\theoremstyle{plain}
\newtheorem{theorem}{Theorem}[section]
\newtheorem{lemma}[theorem]{Lemma}

\newtheorem{proposition}[theorem]{Proposition}
\newtheorem{corollary}[theorem]{Corollary}

\theoremstyle{definition}

\newtheorem{definition}[theorem]{Definition}

\theoremstyle{remark}
\newtheorem{remark}[theorem]{Remark}
\newtheorem{example}[theorem]{Example}

\numberwithin{equation}{section}

\begin{document}

\title[A spectral sequence for groups and groupoids]{A
  spectral sequence to compute $L^2$-Betti numbers of groups and groupoids}  

\thanks{R.S. acknowledges support by DFG grant SA  1661/1-1 and thanks
  the Max-Planck-Institute in Bonn for its 
  hospitality during the initial stage of this project. A.T. thanks the Graduiertenkolleg "Gruppen und Geometrie", G\"ottingen}

\author{Roman Sauer}
\address{University of Chicago, Chicago, USA}
\curraddr{Mathematisches Institut der WWU M\"unster, Einsteinstr.~62, 48149 M\"unster, Germany}
\email{sauerr@uni-muenster.de}
\urladdr{www.romansauer.de}

\author{Andreas Thom}
\address{Mathematisches Institut der Georg-August Universit\"at G\"ottingen, Bunsenstr. 3-5, 37073 G\"ottingen, Germany}
\email{thom@uni-math.gwdg.de}
\urladdr{www.uni-math.gwdg.de/thom}


\begin{abstract}
We construct a spectral sequence for 
$L^2$-type cohomology groups of discrete measured groupoids. 
Based on the spectral sequence, we prove the Hopf-Singer conjecture 
for aspherical manifolds with poly-surface fundamental groups. More generally, we obtain 
a permanence result for the Hopf-Singer conjecture under taking fiber bundles whose base 
space is an aspherical manifold with poly-surface fundamental group. 
As further sample applications of the spectral sequence, 
we obtain new vanishing theorems   
and explicit computations of $L^2$-Betti numbers of 
groups and manifolds and obstructions to the existence of 
normal subrelations in measured equivalence relations. 
\end{abstract}
\maketitle

\section{Introduction and Statement of Results}\label{sec:introduction}

The aim of this work is to construct, starting with a short 
exact sequence (later called \emph{strong extension})
of discrete measured groupoids, a spectral sequence 
for $L^2$-type cohomology groups.
For this, we are using a blend of tools from 
homological algebra and ergodic theory.

\smallskip

Gaboriau introduced and studied 
the notion of $L^2$-Betti numbers of measured equivalence relations~\cite{gaboriau(2002b)}, which proved to be very fruitful, 
especially in applications to von Neumann algebras via the work of Popa. 
Subsequently, more algebraic definitions were 
developed~\cites{sauer(2005),thoml2} 
that build upon L\"uck's algebraic theory of 
$L^2$-Betti numbers~\cite{lueck(2002)}. We work with the latter since they are especially suited 
for our purposes. 

\smallskip

Gaboriau's and L\"uck's works 
add quite different computational 
techniques to the theory of $L^2$-Betti numbers: For instance, 
Gaboriau's theory allows to exploit the fact that the ergodic dimension 
of a group might be much smaller than its cohomological dimension; L\"uck's 
algebraic theory, on the other hand, allows to use the power of standard 
homological algebra in computations of $L^2$-Betti numbers (another 
algebraic $L^2$-theory is due to Farber~\cite{farber}). 

The motivation for our spectral sequence 
was to combine these computational advantages. The reader may wonder  
whether the generality of the language of 
groupoids is necessary; we will present 
computations (Corollaries~\ref{cor:group vanishing in degree plus one}, 
\ref{cor:computation of extensions of free products}, \ref{cor:app to hopf-singer}) for groups and manifolds that, in their 
proofs, make use of measured equivalence relations. The class of measured 
equivalence relations, however, is not closed under taking quotients, 
unlike the class of discrete measured groupoids; so it turns out that 
it is necessary and most natural to work with groupoids in our context. 

\smallskip

We refer the reader for a more detailed 
review of used tools, concepts, and notation to Section~\ref{sec:review}. 

\smallskip

A central notion is that of a \emph{strongly normal subgroupoid} of a
discrete measured groupoid. Building 
on~\cite{feldman+sutherland+zimmer}, 
we see in Section~\ref{sec:quotients and normality} that a
strongly normal subgroupoid (Definition~\ref{def:normality}) allows a
quotient construction, and every strongly normal subgroupoid appears
as the kernel of a quotient map. One calls an ergodic $\calg$ a \emph{strong
  extension} of $\cals$ and $\calq$ if $\cals\subset\calg$ is a
strongly normal subgroupoid and $\calq$ the corresponding quotient
groupoid. Several examples of strong extensions are discussed in
Section~\ref{sec:examples}.

For a discrete measured groupoid $\calg$, we define cohomology groups
$H^\ast(\calg, M)$ with coefficients in a module $M$ over the groupoid
ring of $\calg$ (see Section~\ref{sec: construction}).  If $M$ is the
algebra of affiliated operators $\calu(\calg)$ of the von Neumann
algebra of $\calg$, then L\"uck's dimension of the
$\calu(\calg)$-module $H^n(\calg,\calu(\calg))$ defines the $n$-th
$L^2$-Betti number $\betti_n(\calg)$ of $\calg$. This definition of
$\betti_n(\calg)$ coincides with the definition given by the first
author in~\cite{sauer(2005)} which itself is an algebraic
formulation of Gaboriau's definition~\cite{gaboriau(2002b)}. 
Although it gives the same $L^2$-Betti numbers, our definition 
of $H^\ast(\calg,M)$ is more involved: we have to take, in some sense, 
a huge model of $H^\ast(\calg,M)$ to overcome some serious, technical difficulties. 

Theorem~\ref{thm:spectral sequence} of Section~\ref{sec: construction} 
is our main result: a Grothendieck spectral sequence that
computes the cohomology of a strong extension of $\cals$ and $\calq$
in terms of the cohomology groups of $\cals$ and $\calq$.

\smallskip

We now present applications of the spectral sequence. The proofs of
the theorems below are found in Section~\ref{sec:applications}.

\subsection{Applications to ergodic theory}
The following theorem generalizes a corresponding result by 
Gaboriau~\cite{gaboriau(2002b)}*{Th\'eor\`eme~6.6} for extension of
groups where the quotient group is amenable. 

\begin{theorem}\label{thm:vanishing with amenable quotient}
Let $$1\rightarrow\cals\longrightarrow\calg\longrightarrow\calq\rightarrow 1$$ 
be a strong 
extension of discrete measured groupoids where $\calq$ is an infinite amenable  discrete measured groupoid. If $\betti_n(\cals)<\infty$,  then
$\betti_n(\calg)=0$. 
\end{theorem}

\begin{corollary}\label{cor:wild subrelations}
Let $\Gamma$ be a countable group with $\betti_1(\Gamma)\ne 0 $. 
Then there is a free, ergodic $\Gamma$-probability space such that 
the associated orbit equivalence relation $\calr$ has a strongly normal 
subrelation $\cals\subset\calr$ of infinite index with $\betti_1(\cals)=\infty$. 
\end{corollary}

\begin{proof}[Proof of Corollary]
Because of $\betti_1(\Gamma)\ne 0$ the group $\Gamma$ does not have 
property $(T)$. By~\cite{schmidt}*{Theorem~1.5} there is an ergodic 
$\Gamma$-probability space $(X,\mu)$ that is not \emph{strongly
  ergodic} (defined in~\cite{schmidt}*{Section~1}). 
By~\cite{schmidt}*{Theorems~2.1 and~2.3} (the details of the proof are 
in~\cite{jones-schmidt}), there is an amenable  
equivalence relation~$\calr_{\rm hyp}$ and a strong surjection 
$\theta:X\rtimes\Gamma\rightarrow\calr_{\rm hyp}$. 
The claim now follows for $\cals = \ker(\theta)$ from 
Theorem~\ref{thm:vanishing with amenable quotient} and 
$\betti_p(X\rtimes\Gamma)=\betti_p(\Gamma)$ 
(see Section~\ref{subsec:identification}). 
\end{proof}

\begin{theorem}\label{thm:betti-less normal subrelations} 
If an ergodic discrete measured groupoid $\calg$ possesses a strongly normal 
subgroupoid $\cals$ such that $\betti_n(\cals)=0$ for every $n \leq d$,
then $\betti_n(\calg)=0$ for every $n \leq d$. 
\end{theorem}

The question of (non-)existence of 
subrelations that are (strongly) normal or split off as a
factor is investigated in several works. For example, in~\cite{feldman+sutherland+zimmer}*{Theorem~4.4}
it is proved that orbit equivalence relations of free II$_1$-actions
of a lattice in a simple, connected, non-compact Lie group never have
infinite, strongly normal, amenable
subrelations. In~\cite{adams}*{Corollary~6.2} it is shown that the
orbit equivalence relation of a type II$_1$-action of a non-elementary
word-hyperbolic group is irreducible, that is, it 
never splits off an infinite, discrete measured equivalence relation as a factor. 

\begin{corollary}\label{cor:amenable groupoid are betti-less}
Let $\calg$ be an infinite, amenable, ergodic discrete measured groupoid. Then
$\betti_p(\calg)=0$ for every $p\ge 0$. 
\end{corollary}

\begin{proof}[Proof of Corollary]
Consider the principal extension $1\rightarrow\calg_{\rm
  stab}\rightarrow\calg\rightarrow\calg_{\rm rel}\rightarrow 1$ of 
Subsection~\ref{subsec:principal extension}. Assume first 
that $\calg_{\rm stab}$ is infinite. Since $\calg$ is amenable, the
stabilizer groups $\calg_x=\{\gamma\in\calg;~r(\gamma)=s(\gamma)=x\}$
are (infinite) amenable 
for a.e.~$x\in\calg^0$. By 
Lemma~\ref{lem: betti number of groupoid with only isotropy}, we have 
$\betti_p(\calg_{\rm stab})=0$ for every $p\ge 0$. By
Theorem~\ref{thm:betti-less normal subrelations} the assertion
follows. 

Suppose now that $\calg_{\rm stab}$ is finite. Then $\calg_{\rm rel}$
has to be infinite since $\calg$ is so. Hence $\calg_{\rm rel}$ is an 
infinite 
hyperfinite equivalence relation. Now the assertion is implied by 
Lemma~\ref{lem: betti number of groupoid with only isotropy} and 
Theorem~\ref{thm:vanishing with amenable quotient}. 
\end{proof}

\begin{theorem}\label{thm:vanishing in degree plus one}
Let 
\[1\rightarrow\cals\longrightarrow\calg\longrightarrow
\calq\rightarrow 1\] 
be a strong
extension of discrete measured groupoids where $\calq$ is an infinite
measured groupoid. Assume that $\betti_p(\cals)=0$ for all $0\le p\le
d$ and $\betti_{d+1}(\cals)<\infty$. Then $\betti_p(\calg)=0$ for all $0\le
p\le d+1$. 
\end{theorem}

As a consequence, we obtain an alternative proof 
of~\cite{bergeron+gaboriau}*{Th\'eor\`eme~5.4}: 

\begin{corollary}[Bergeron-Gaboriau]\label{cor:first betti and stabilizers}
Let $\Gamma$ be a countable group with $\betti_1(\Gamma)\ne 0$. Let
$(X,\mu)$ be an ergodic $\Gamma$-probability space. Then either 
the stabilizer $\Gamma_x$ is finite for $\mu$-a.e.~$x\in X$, or
$\betti_1(\Gamma_x)=\infty$, thus $\Gamma_x$ is not finitely
generated, for $\mu$-a.e.~$x\in X$. 
\end{corollary}

\begin{proof}[Proof of Corollary]
By ergodicity, either $\Gamma_x$ is finite almost everywhere or infinite 
almost everywhere. The function $x\mapsto\betti_1(\Gamma_x)$ is 
measurable~\cite{bergeron+gaboriau}*{Th\'eor\`eme~5.7}. Thus, by ergodicity again, 
either $\betti_1(\Gamma_x)=\infty$ almost everywhere or 
$\betti_1(\Gamma_x)<\infty$ almost everywhere. 
Apply Theorem~\ref{thm:vanishing in degree plus one} 
to the principal extension
(Subsection~\ref{subsec:principal extension}) in combination with
Lemma~\ref{lem: betti number of groupoid with only isotropy}. 
Note that the first $L^2$-Betti number of a finitely generated group is finite. 
\end{proof}

\begin{remark}[Ergodicity hypothesis]\label{rem:remark on ergodicity}
We remark that 
the notion of \emph{strong extension} is based upon 
ergodicity. Nevertheless, it is possible to drop the  
assumption on ergodicity in 
several of the results above. In fact, if $\int_X\mu_xd\mu(x)=\mu$ 
is the ergodic decomposition of $(\calg,\mu)$, then it is possible 
to show that 
\[\betti_p(\calg,\mu)=\int_X\betti_p(\calg_x,\mu_x)d\mu(x).\]
\end{remark}

\subsection{Applications to $L^2$-Betti numbers of groups and manifolds}

L\"uck~\cite{lueck(2002)}*{Theorem~7.2 on p.~294} 
proved the following corollary under the additional assumption 
that the quotient group has a non-torsion element or finite subgroups of 
arbitrarily high order. Then, using equivalence relations, 
Gaboriau~\cite{gaboriau(2002b)}*{Th\'eor\`eme~6.8}
gave a proof in degree~$1$, i.e.~for the first $L^2$-Betti number, 
only assuming the quotient to be infinite. We also 
mention that~\cite{bourdon}*{Corollary~1} 
gives a very elementary proof for the degree~$1$ case if the quotient
has a non-torsion element. It is remarked therein that it is \emph{a
challenging, and vaguely irritating question to find a purely
cohomological proof of Gaboriau's result}. Up to now, there is no
proof that does not use measured equivalence relations. 
The following corollary of 
Theorem~\ref{thm:vanishing in degree plus one} 
generalizes the aforementioned results to all 
degrees without further assumptions on the quotient. 

\begin{corollary}\label{cor:group vanishing in degree plus one}
Let $\Lambda\subset\Gamma$ be a normal subgroup of infinite index. Suppose 
that $\betti_p(\Lambda)=0$ for $0\le p\le d-1$ and
$\betti_d(\Lambda)<\infty$. Then $\betti_p(\Gamma)=0$ for $0\le p\le
d$. 
\end{corollary}

\begin{proof}[Proof of Corollary]
One can find free, ergodic measure-preserving actions of
$\Lambda,\Gamma$ and $Q=\Gamma/\Lambda$ on probability spaces such
that the associated orbit equivalence relations form a strong
extension (see Subsection~\ref{subsec:example with groups by gaboriau}). 
Then apply Theorem~\ref{thm:vanishing in degree plus one} 
and~\eqref{eq:betti identification} in Section~\ref{sec: construction}. 
\end{proof}

The notion of \emph{measure equivalence} was introduced by 
Gromov and, for the first time, 
gained prominence in the work of 
Furman~\cite{furman(1999a)}*{Definition~1.1}:

\begin{definition}\label{def:measure equivalence}
  Two countable groups $\Gamma$ and $\Lambda$ are called \emph{measure equivalent} 
if there exists a non-trivial measure space $(\Omega,\mu)$ on which $\Gamma\times \Lambda$ acts such that the
  restricted actions of $\Gamma=\Gamma\times\{1\}$ and $\Lambda=\{1\}\times\Lambda$ have
  measurable fundamental domains $X\subset\Omega$ and
  $Y\subset\Omega$ with $\mu(X)<\infty$ and $\mu(Y)<\infty$. 
   The space $(\Omega,\mu)$ is called a
  \emph{measure coupling} between $\Gamma$ and $\Lambda$.
\end{definition}

We denote by $\cd_\bbC(\Gamma)$ the 
\emph{cohomological dimension} of a group $\Gamma$ over $\bbC$, 
i.e.~the projective dimension of $\bbC$ as a $\bbC \Gamma$-module. 
The following theorem is a consequence of the more general 
Theorem~\ref{secondvariant} and Lemma~\ref{newlemma}. 

\begin{theorem}\label{thm:vanishing by small ergodic dimension}
Let $1\rightarrow\Lambda\rightarrow\Gamma\rightarrow Q_0\rightarrow 1$
be a short exact sequence of groups. 
Suppose 
that $\betti_p(\Lambda)=0$ for $p>m$. 
Let $Q_1$ be a group that is measure equivalent to $Q_0$. 
Let $n=\cd_\bbC(Q_1)$. Then 
$\betti_p(\Gamma)=0$ for $p>n+m$. 
\end{theorem}

\begin{remark} \label{remmcd}
In order to capture the strength of the method of proof that leads to the above theorem, 
we introduce the concept of \emph{measurable cohomological dimension} in Section~\ref{singer}. It is then clear that the only relevant assumption is that the measurable cohomological dimension of the quotient group is bounded by $n$. We will prove this fact in Theorem~\ref{secondvariant}.
\end{remark}

\begin{remark}
A typical situation where the quotient group is measure equivalent 
to a group with a lower cohomological dimension is the following: Let 
$\Gamma$ be a cocompact lattice in $\SL(n,\bbR)$. Then $\SL(n,\bbR)$ endowed with the 
left multiplication action of $\Gamma$ and the right multiplication action of 
$\SL(n,\bbZ)$ is a measure coupling with respect to the Haar measure. 
The rational cohomological 
dimension of $\Gamma$ equals the dimension of the 
associated symmetric space $\SL(n,\bbR)/\SO(n,\bbR)$; but the rational cohomological dimension 
of $\SL(n,\bbZ)$ is 
\[\cd_\bbQ\bigl(\SL(n,\bbZ)\bigr)=\dim\bigl(\SL(n,\bbR)/\SO(n,\bbR)\bigr)-(n-1)\]
by a result of Borel-Serre~\cite{borel-serre}. 
\end{remark} 

We present the following sample 
application of Theorem~\ref{thm:vanishing by small ergodic dimension}, for which we do 
not know an alternative proof that does not use measured groupoids. 

\begin{corollary}\label{cor:computation of extensions of free products}
Let $A_1,\ldots,A_k$ and $B_1,\ldots,B_l$ be infinite amenable
groups. Let $\Gamma$ be an extension of the type 
\begin{equation*}
1\rightarrow A_1\ast\cdots\ast A_k\rightarrow\Gamma\rightarrow
B_1\ast\cdots\ast B_l\rightarrow 1. 
\end{equation*}
Then 
\begin{equation*}
\betti_p(\Gamma)=\begin{cases}
                  (k-1)(l-1) &\text{ if $p=2$,}\\
                  0 &\text{ otherwise.}
                 \end{cases}
\end{equation*}
\end{corollary}

\begin{proof}[Proof of Corollary]
Since $\betti_p(A_i)=0$ for every $p\ge 0$ and $i\in\{1,\ldots,k\}$ 
by~\cite{lueck(2002)}*{Theorem~6.37 on p.~259}, 
we obtain by the Mayer-Vietoris sequence for $l^2$-homology 
that $\betti_p(A_1\ast\cdots\ast A_k)=0$ whenever $p>1$. 

Since $B_1\ast\cdots\ast B_l$ and the free group $F_l$ of rank $l$ 
are measure equivalent~\cite{gaboriau-examples}*{Section~2.2}, 
it follows that $\betti_p(\Gamma)=0$ for $p>2$ by the
previous theorem. Since $\Gamma$ and $A_1\ast\cdots\ast A_k$ 
are infinite, we have
$\betti_0(\Gamma)=\betti_0(A_1\ast\cdots\ast A_k)=0$, thus,   
by Corollary~\ref{cor:group vanishing in degree plus one}, 
$\betti_1(\Gamma)=0$. Since the Euler characteristic $\chi$ is
multiplicative for extensions and 
$\chi(\Gamma)=\sum_{p\ge 0}(-1)^p\betti_p(\Gamma)$, we obtain that 
\begin{equation*}
\betti_2(\Gamma)=\chi(\Gamma)=\chi(A_1\ast\cdots\ast
A_k)\chi(B_1\ast\cdots\ast B_l)= (k-1)(l-1).\qedhere
\end{equation*}
\end{proof}

Next we turn our attention to a central conjecture for $L^2$-Betti numbers, 
the \emph{Hopf-Singer Conjecture}; it 
predicts that for a closed aspherical manifold $M$ we have 
$\betti_p(\widetilde{M})=0$ provided $2p\ne \dim(M)$. 
For a survey of known results we refer the reader
to~\cite{lueck(2002)}*{Chapter~11}. 

A group is said to be a \emph{poly-surface group}, if it has a series of normal subgroups such that the subquotients are fundamental groups of closed oriented surfaces (see Definition~\ref{defpolysurf} for more details). The cohomological dimension of such a group is precisely twice the length of the series and the Euler characteristic is the product of the individual Euler characteristics of the subquotients. 

In Section~\ref{singer}, we study the measurable cohomological dimension (which is closely related to Gaboriau's ergodic dimension~\cite{gaboriau(2002b)}) in more detail and show, that the measurable cohomological dimension of a poly-surface group is the length of its defining series of normal subgroups, i.e.\ only half of the expected number.

The following theorem is obtained as Corollary~\ref{polysinger} in 
Subsection~\ref{sub:singer condition}. 

\begin{theorem}
	The Hopf-Singer conjecture holds for any closed aspherical manifold with 
	poly-surface fundamental group. 
\end{theorem}

We also obtain a permanence result for the Hopf-Singer conjecture: 

\begin{theorem}\label{cor:app to hopf-singer}
Let $M$ be a closed, aspherical, $2n$-dimensional manifold that
satisfies the Hopf-Singer Conjecture, that is, $\betti_p(\widetilde{M})=0$
unless $p=n$. Let $L$ be a closed orientable aspherical manifold of dimension $2m$ with a poly-surface fundamental group. 
If $N$ is the total space of an orientable fiber bundle over $L$ with fiber $M$, then 
\begin{equation*}
\betti_p(\widetilde{N})=\begin{cases}
               \betti_n(\widetilde{M})\cdot \chi(\pi_1(L))&\text{
                 if $p= m+n$,}\\
               0 &\text{ otherwise.}
            \end{cases}
\end{equation*}
In particular, $N$ satisfies the Hopf-Singer Conjecture, too. 
\end{theorem}

\begin{proof}
The manifold $N$ is closed, orientable and aspherical. Let $\Gamma=\pi_1(N)$,
and $\Lambda=\pi_1(M)$. From the fiber bundle we get an extension of
groups 
\begin{equation*}
1\rightarrow\Lambda\rightarrow\Gamma\rightarrow\pi_1(L)\rightarrow 1. 
\end{equation*}
Now, $\pi_1(L)$ satisfies Singer's condition (see Definition \ref{defsinger}) by Theorem \ref{polysinger}. Hence, the group $\pi_1(L)$ has measurable cohomological dimension $m$. 
By Theorem~\ref{secondvariant}, 
$\betti_p(\widetilde{N})=\betti_p(\Gamma)=0$ for $p>m+n$. 
Thus, by Poincar\'e duality, $\betti_p(\widetilde{N})=0$ unless $p=m+n$. 
This yields $\betti_{m+n}(\widetilde{N})=\pm\chi(N)$, and from the
multiplicativity of $\chi$ for fiber bundles the claim follows. 
\end{proof}

\section{Overview of used concepts and tools}\label{sec:review}

\subsection{Standard Borel and measure spaces} 
All measurable spaces in this work are 
\emph{standard Borel spaces} unless stated otherwise. 
Maps between standard Borel spaces are measurable unless stated
otherwise. 
Our background references for standard Borel spaces
are~\cites{cohn,kechris}; we recall here some 
basic notions and facts (see also~\cite{glasner}*{p.~51/52}). 

\smallskip

We use the terms \emph{measurable} and
\emph{Borel} interchangeably for maps between or subsets of standard Borel
spaces. A measure on a standard Borel space is called \emph{Borel
  measure}. A \emph{partition} of standard 
Borel space $X$ is a countable family
$(X_i)_{i\in I}$ of pairwise disjoint Borel subsets such that
$X=\bigcup_{i\in\bbN} X_i$. 
A \emph{Borel isomorphism} $f:X\rightarrow Y$ between standard
Borel spaces is a bijective Borel map. 
Inverses of Borel isomorphisms are Borel, and 
Borel subsets of a standard Borel space are again standard Borel. 

\smallskip

The following result
is a fundamental tool to which we refer throughout the paper as the
\emph{selection theorem} (see~\cite{kechris}*{theorem~18.10 on p.~123}). 
\begin{theorem}[Selection Theorem]
Let $f:X\rightarrow Y$ be a Borel map between standard Borel spaces
whose fibers are countable. Then $f(X)\subset Y$ is Borel, and
there is a partition $X=\bigcup_{i\in\bbN}X_i$ such that each
$f\vert_{X_i}$ is injective. 
\end{theorem}
A \emph{measure space} $(X,\mu)$ is by definition a standard Borel
space $X$ equipped with a Borel measure $\mu$. A \emph{probability
  space} is a measure space whose measure is a probability measure, that
is, has total measure $1$. 
A \emph{measure isomorphism} $f:(X,\mu_X)\rightarrow (Y,\mu_Y)$ is a
measure-preserving Borel
map with the property that there are Borel subsets $A\subset X$ and $B\subset Y$ with 
$\mu_X(X-A)=\mu_Y(Y-B)=0$ such that $f\vert_A$ is a Borel isomorphism
$A\rightarrow B$. If $(X,\mu)$ is continuous, that is, $\mu(\{x\})=0$
for every $x\in X$, then there is Borel isomorphism $f:X\rightarrow
[0,1]$ with $f_\ast\mu=\mu\circ f^{-1}=\lambda$. Here, $\lambda$ denotes the Lebesgue measure. 

\smallskip

Another important tool~\cite{glasner}*{Theorem~3.22 on p.~72} is 
\begin{theorem}[Measure disintegration]
Let $(X,\mu)$ and $(Y,\nu)$ be probability spaces and
$\pi:X\rightarrow Y$ a Borel map such that $\pi_\ast\mu=\nu$. 
Then there is a map $y\mapsto\mu_y$ that associates to every $y\in Y$
a probability measure $\mu_y$ on $X$ such that 
\begin{enumerate}[i)]
\item For every Borel subset $A\subset X$, the function
$y\mapsto \mu_y(A)$ is Borel. 
\item For $\nu$-a.e. $y\in Y$, $\mu_y(\pi^{-1}(y))=1$. 
\item $$\mu=\int_Y\mu_yd\nu(y).$$
\end{enumerate}
\end{theorem}

\subsection{Discrete measured groupoids}
The standard reference for measurable groupoids is~\cite{anan}. 
The \emph{source and range maps} of a groupoid are denoted by $s$ and $r$,
respectively. We use superscript $0$ to denote 
the \emph{unit space} of a groupoid, as in $\calg^0$ for the groupoid
$\calg$. Our convention is that a \emph{subgroupoid} of a groupoid 
has the same unit space. 

\smallskip

A \emph{discrete measurable groupoid} $\calg$ is a groupoid $\calg$ 
equipped with the structure of a standard Borel space such that 
$\calg^0\subset\calg$ is a Borel subset, all the structure maps are
Borel, and $s^{-1}(\{x\})$ is countable for all $x\in\calg^0$. 
Let $c^s_x$ and $c^r_x$ denote the counting measures on $s^{-1}(x)$
and $r^{-1}(x)$, respectively. 

\smallskip

A \emph{discrete measured groupoid} $(\calg,\mu)$ 
is a discrete measurable groupoid $\calg$ together with a Borel
measure $\mu$ on $\calg^0$ such that 
\begin{equation*}
\int_{\calg^0}c^s_xd\mu(x)=\int_{\calg^0}c^r_xd\mu(x)
\end{equation*}
as measures on $\calg$. This measure on $\calg$, 
which extends the one on $\calg^0$, is also denoted by
$\mu$. 

\smallskip

Moreover, $(\calg,\mu)$ is called \emph{infinite} if 
$s^{-1}(x)$ is infinite for all $x\in\calg^0$ in a subset of 
positive measure. 

\smallskip

The discrete measured groupoid 
$(\calg,\mu)$ is said to be \emph{ergodic} if one (thus, all) of the following 
equivalent conditions hold: 
\begin{enumerate}[i)]
\item Any function $f:\calg^0\rightarrow\bbR$ that is
  \emph{$\calg$-invariant} (that is, $s\circ f=r\circ f$) is $\mu$-a.e.\ constant. 
\item For any Borel subset $A\subset\calg^0$ of positive measure, the so-called 
\emph{saturation} 
$A^\calg\defq\{x\in\calg^0;~\exists\gamma\in\calg: s(\gamma)\in A, r(\gamma)=x\}$ 
has full measure.  
\item For any two Borel subsets $A,B\subset\calg^0$ of positive
  measure there exists a Borel subset $E\subset\calg$ such that
  $s(E)\subset A$, $r(E)\subset B$ and $\mu(s(E))>0$, $\mu(r(E))>0$. 
\end{enumerate}
Each discrete measured groupoid $(\calg,\mu)$ has an \emph{ergodic decomposition}, 
that is, there is an disintegration map $(\calg^0,\mu)\rightarrow (X,\nu)$ 
such that $(\calg,\mu_x)$ is 
for $\nu$-a.e.~$x\in X$ ergodic~\cite{hahn}*{Theorem~6.1}. 
The measure space $(X,\nu)$ 
is called the \emph{space of ergodic components}. 

\smallskip

A Borel map $\phi:(\calg, \mu)\rightarrow (\calg',\nu)$ is called 
a \emph{homomorphism} if $\phi$ is a map of groupoids and
$f_\ast\mu=\nu$. 

\smallskip

Let $(X,\mu)$ be a probability space with a probability measure
$\mu$, and let $\Gamma\times X\rightarrow X$ be a measure preserving
(m.p.)\ group action. We denote by $(X\rtimes \Gamma, \mu)$ the
translation groupoid, i.e.~the groupoid with total space
$X\times\Gamma$, base space $X$ and where $s=\pi_X$ and
$r:X\rtimes\Gamma\rightarrow X$ is defined to be the action of $\Gamma$ on $X$. 

For the definition of an \emph{amenable} or, equivalently, \emph{hyperfinite} discrete 
measured groupoid we refer to~\cite{tak3}*{Chapter~XIII, \S 3}. 
\subsection{Spectral sequences}

Weibel's book \cite{weibel} is a standard reference for all the 
homological algebra that we need. We restate
\cite{weibel}*{Theorem~5.8.3 on p.~158} for 
the convenience of the reader.

\begin{theorem}[Grothendieck] \label{grothendieck} Let $\cala,\calb$
  and $\calc$ be abelian categories, such that both $\cala$ and
  $\calb$ have enough injective objects. Let $G \colon \cala \to
  \calb$ and $F \colon \calb \to \calc$ be left exact functors, such
  that that $G$ sends injectives to injectives. Then, given $A \in
  \cala$, there exists a first quadrant spectral sequence:
$$E_2^{pq} = (R^pF)(R^qG)(A) \Longrightarrow R^{p+q}(FG)(A).$$
\end{theorem}

Here $R^pF$, for $p \geq 0$, denotes the $p$-th right derived functor
of $F$.  This very general spectral sequence can be used to construct
the classical Hochschild-Serre spectral sequence in group cohomology. 
The core of our work consists in constructing suitable functors $F$ and $G$ and 
verifying the above hypothesis in the setting of discrete 
measured groupoids. 

In view of the above assumptions, it is useful to have 
a criterion that ensures that a functor preserves injective objects: 

\begin{lemma}[\cite{weibel}*{Theorem~2.6.1 on p.~50}] \label{inj-criterion}
Let $G\colon \cala \to \calb$ be a right exact functor between abelian categories. If $G$ has an exact left-adjoint functor,
then it sends injectives to injectives.
\end{lemma}

\subsection{$L^2$-Betti numbers} The standard reference 
for $L^2$-Betti numbers is L\"uck's book~\cite{lueck(2002)}. 
In this paper we only refer to the homological-algebra-type definition 
of $L^2$-Betti numbers in the context of 
groups~\cite{lueck(2002)}*{Chapter~6} and 
discrete measured groupoids~\cite{sauer(2005)} (see also~\cite{thomrank}). 

\smallskip

Let $(M,\tau)$ be a finite von Neumann algebra with a fixed trace
$\tau\colon M \to \bbC$. 
L\"uck introduced a dimension function $L \mapsto \dim_{(M,\tau)}(L)$ for an arbitrary 
$M$-module  
$L$ with very nice properties like additivity for arbitrary modules. 
If the context is clear, 
we may also omit the trace $\tau$ or the whole subscript 
in $\dim_{(M,\tau)}$. 

\smallskip

If $M=L(\Gamma)$ is the group von Neumann algebra of a group $\Gamma$ 
with its standard trace, then the $p$-th $L^2$-Betti number of $\Gamma$ is 
defined as 
\[\betti_p(\Gamma)\defq\dim_{L(\Gamma)}\Tor_p^{\bbC\Gamma}\bigl(\bbC,L(\Gamma)\bigr)
   =\dim_{L(\Gamma)}H_p\bigl(\Gamma, L(\Gamma)\bigr)\in [0,\infty].\]
If $\calg$ is a discrete measured groupoid, then one defines 
in a similar way 
\[\betti_p(\calg)\defq
\dim_{L(\calg)}\Tor_p^{\calr(\calg)}\bigl(L^\infty(\calg^0),L(\calg)\bigr)
\in [0,\infty].\]
For the definition of $\calr(\calg)$ and $L(\calg)$ see 
Section~\ref{subsec:modules}.

\subsection{Operators affiliated with a finite von Neumann algebra}

Let $(M,\tau)$ again be a finite von Neumann algebra with a fixed trace
$\tau$. Denote by $L^2(M,\tau)$ the GNS-construction
with respect to $\tau$. There is a natural algebra
$\calu(M,\tau)$ of closable, densely defined and unbounded operators
on $L^2(M,\tau)$, which are affiliated with the algebra $M$. For
details, we refer to \cite{tak2}*{Chapter IX}. 

\smallskip

There is  
also the notion of dimension 
$\dim_{\calu(M,\tau)}(M)\in [0,\infty]$ for an arbitrary 
$\calu(M,\tau)$-module $M$~\citelist{\cite{reich}\cite{lueck(2002)}*{Chapter~8}}. As it turns out,
a useful dimension function is obtained by restricting the module structure to $M$ and taking the dimension
of the $\calu(M,\tau)$-module as an $M$-module.
Again, we may omit $\tau$ or the whole subscript in $\dim_{\calu(M,\tau)}$ if the context is clear. 

\smallskip

This algebra has been studied in connection with $L^2$-invariants in
\cites{reich, thoml2}. 
It shares very nice ring-theoretic
properties, which we want to summarize in the sequel:
\begin{enumerate}[i)]
\item $\calu(M,\tau)$ is \emph{von Neumann regular}, i.e.\ all modules are
  flat,
\item $\calu(M,\tau)$ is left and right \emph{self-injective},
  i.e. $\calu(M,\tau)$ is injective as left and as right module over
  itself,
\item $\iota \colon M \to \calu(M,\tau)$ is a flat ring extension,
  i.e.\ $L \mapsto \calu(M,\tau) \otimes_M L$ is exact,
\item $\dim_{\calu(M,\tau)} \calu(M,\tau) \otimes_M L = \dim_{(M,\tau)}L$, for every
  $M$-module $L$,
\item $\dim_{\calu(M,\tau)}\hom_{\calu(M,\tau)}(L,\calu(M,\tau)) = 
      \dim_{\calu(M,\tau)} L$, for every $\calu(M,\tau)$-module $L$.
\item If $\dim_{\calu(M,\tau)} L =0$, then $\hom_{\calu(M,\tau)}(L,\calu(M,\tau)) =0$. 
\end{enumerate}
For proofs, we refer to \cites{reich,thoml2} and the references
therein.

\subsection{Localization and completion}
Let $(M,\tau)$ be as above. Let $\cala$ be an abelian category with a
faithful functor $F\colon \cala \to \modules{M}$. Assume that the
functor $F$ preserves limits and co-limits. In the category
of $M$-modules, the sub-category of zero-dimensional modules is a
Serre sub-category. Moreover, the full sub-category $\cals \subset \cala$,
given by those modules, that map to zero-dimensional modules, forms a
Serre sub-category as well.

\begin{example}
  The example to have in mind, is the category of $L^{\infty}(X)
  \rtimes \Gamma$-modules with its forgetful functor to
  $L^{\infty}(X)$-modules.
\end{example}

Some of our computations will be carried
out in the \emph{localized}
category $\cala/\cals$, which we also denote by $\cala_{\rm
  loc}$. This category is abelian and its properties, for the examples
of special interest, were studied in~\cite{thomrank}, where the second
author showed that it has enough projective objects and naturally
embeds in $\cala$ as the sub-category of those modules, which are
\emph{complete} with respect to the rank metric, see~\cite{thomrank} for
details. We recall some of the results in Section~\ref{sec: construction}.

\section{Quotients and Normality}\label{sec:quotients and normality}

In this section, we provide the technical underpinning of the concept
of (strong) extension of discrete measurable
groupoids. 

\subsection{Ergodic discrete measurable groupoids}

\begin{lemma}\label{lem:intersection large}
Let $(\calg,\mu)$ be an ergodic discrete measured groupoid with  
atom-free $\mu$, that is, 
$\mu(\{x\})=0$ for every $x\in\calg^0$. 
Let $A, B\subset \calg^0$ be Borel subsets of positive measure. Then 
\begin{equation*}
\mu\bigl(s^{-1}(A)\cap r^{-1}(B)\bigr)=\infty.
\end{equation*}
\end{lemma}

\begin{proof}
We start with a general observation about ergodic groupoids: 
The function $\calg^0\rightarrow\bbZ\cup\{\infty\}, 
x\mapsto\#r^{-1}(x)=\#s^{-1}(x)$ is a.e.\ constant because of
ergodicity. If $\calg^0$ is atom-free, then $\# s^{-1}(x)=\infty$ 
for a.e.~$x\in \calg^0$ since otherwise there existed $n\ge 1$ 
and $A\subset \calg^0$
with $\mu(A)<1/n$ and $\#s^{-1}(x)\le n$ for $x\in A$, implying 
$\mu(A^\calg)<1$. 

Again by ergodicity, 
we can pick a set $E\subset\calg$ with $A'\defq s(E)\subset A$ and 
$B'\defq r(E)\subset B$ such that 
$\mu(A'),\mu(B')>0$. 
By the selection 
theorem we can assume that there is such $E$ with 
injective $s\vert_E$ and $r\vert_E$.  
Denote by $f:B'\rightarrow E$ the inverse of $r:E\rightarrow B'$. 
Consider the Borel map 
\begin{equation*}
  \phi:s^{-1}(A')\cap r^{-1}(B')\rightarrow s^{-1}(A')\cap
  r^{-1}(A')=\calg_{A'},
  \phi(\gamma)=f(r(\gamma))^{-1}\circ\gamma.
\end{equation*}
Notice that the restricted groupoid $\calg_{A'}\defq s^{-1}(A')\cap
r^{-1}(A')$ is ergodic if $\calg$ is so.

Since $\phi$ is a map over $A'$ with respect to the source maps and
fiberwise bijective, we obtain, with the general observation above,
that
\begin{equation*}
  \mu\bigl(s^{-1}(A)\cap r^{-1}(B)\bigr)\ge\mu\bigl(s^{-1}(A')\cap
  r^{-1}(B')\bigr)=\mu\bigl(\calg_{A'}\bigr)=\infty. \qedhere
\end{equation*}
\end{proof}

The following lemma is certainly well known 
but we failed to find a reference. 

\begin{lemma}\label{lem:selection theorem with injective range}
Let $(\calg,\mu)$ be an ergodic discrete measured groupoid. 
Then there is a countable set $I$ and 
a measure isomorphism $\phi:\calg^0\times I\rightarrow\calg$ 
such that 
$s\circ\phi=\pr_\calg^0$, where 
$\pr_\calg^0:\calg^0\times I\rightarrow \calg^0$ is the projection 
and $\calg^0\times I$ is endowed with the product of $\mu$ and the 
counting measure on $I$. 
Further, for every $i\in I$, the map 
$\calg^0\rightarrow \calg^0, x\mapsto r(\phi(x,i))$, is a measure 
isomorphism. 
\end{lemma}
\begin{proof}
By ergodicity, 
$(\calg^0,\mu)$ is either discrete (thus $\calg$ is finite) or atom-free. 
We leave the easy proof of the first case to the reader and proceed to
the atom-free case. 
The following auxiliary fact is needed for the proof.\smallskip\\
\textit{Claim:} Let $F\subset\calg$ be a Borel subset of finite
measure. 
Let $E\subset\calg\backslash F$ be a Borel subset on which $r$
and $s$ are injective, then there 
is a Borel subset $D\subset\calg\backslash F$ 
containing $E$ such that $r, s$ are
injective on $D$ and $s(D)=\calg^0$ up to null sets. Note that by
$\calg$-invariance of $\mu$, this also implies that $r(D)=\calg^0$ up to
null sets.
 \smallskip\\
The set $\calb$ of Borel sets $D$ such that $E\subset
D\subset\calg\backslash F$ and
$r\vert_D,s\vert_D$ are injective contains $E$ and is partially
ordered as follows: $D_1\le D_2$ if and only if 
$D_1\subset D_2$ up to null sets. 
Let $\calt\subset\calb$ be a totally ordered subset. To apply
Zorn's lemma later on, we show that $\calt$ has an upper bound in
$\calb$. Let $r=\sup_{D\in\calt}\mu(s(D))\in (0,1]$. Pick a countable 
family $\{D_n\}_{n\in\bbN}$ in $\calt$ with $D_n\le D_{n+1}$ and
$\mu(s(D_n))\rightarrow r$. Let $D=\bigcup D_n$. Upon subtracting a
suitable null
set from $D$ we can ensure that for all $x,y\in D$ there is
$n_0\in\bbN$ with $x,y\in D_{n_0}$. Thus $r\vert_D, s\vert_D$ are
injective, and $D$ is an upper bound of $\calt$. 
By Zorn's lemma the set $\calb$ possesses a maximal element
$D_{max}\in\calb$. Suppose $\mu(s(D_{max}))<1$, thus $\mu(s(D_{max}))<1$. Let
$A\subset \calg^0\backslash s(D_{max})$ and 
$B\subset \calg^0\backslash r(D_{max})$ be
subsets of positive measure. Since $\mu(F)<\infty$ by assumption, 
Lemma~\ref{lem:intersection large} implies that 
$\mu(s^{-1}(A)\cap r^{-1}(B)\cap (\calg\backslash F))=\infty$. 
In particular, there is 
$E\subset\calg$ of positive measure 
such that $s(E)\subset A$ and $r(E)\subset B$. 
Once more using the selection theorem, we can assume that 
there is such $E$ with $r,s$ being injective on $E$. 
Since $E\cup D_{max}$ would contradict
maximality of 
$D_{max}$, it follows that $\mu(s(D_{max}))=\mu(r(D_{max}))=1$. 

\smallskip

Let us continue with the proof of the lemma. 
Using the selection theorem, pick a countable
Borel partition $\{E_n\}_{n\in\bbN}$ of $\calg$ such that $r,s$ are injective on 
each $E_n$. We construct inductively 
Borel subsets $D_n\subset\calg$ such that
\begin{enumerate}[i)]
\item $E_n\subset\bigcup_{i=1}^nD_i$ for every $n\ge 1$, 
and the union is disjoint up to null sets, 
\item $r$ and $s$ are injective on $D_n$ for every $n\ge 1$, 
\item $r(D_n)=s(D_n)=\calg^0$ up to null sets. 
\end{enumerate}
According to the claim above there is such $D_1$, and 
for $n>1$ there is a Borel subset 
$D_n\subset\calg\backslash\bigcup_{i=1}^{n-1}D_i$ 
containing $E_n\backslash\bigcup_{i=1}^{n-1}D_i$ such
that $r,s$ are injective on $D_n$ and $s(D_n)=\calg^0$ up to null sets. 
Since $E_n\subset\bigcup_{i\ge 1}D_i$ for every $n\in\bbN$, we have 
$\calg=\bigcup_{i\ge 1} D_i$. The (inverse of the) 
desired isomorphism is now given as 
\begin{equation*}
\calg\rightarrow \calg^0\times\bbN, \gamma\mapsto (s(\gamma), n)\text{ if
  $\gamma\in D_n$.}\qedhere
\end{equation*}
\end{proof}

\begin{definition}\label{def:endomorphisms}
Let $\cals\subset\calg$ be a subgroupoid of a discrete
measured groupoid $\calg$. 
Let $A\subset\calg^0$ be a Borel subset. 
Let $\emor_\cals(\calg)\vert_A$ denote the set of all Borel maps
$\phi:A\rightarrow\calg$ such that 
\begin{enumerate}[i)]
\item $s(\phi(x))=x$ for a.e.~$x\in A$, 
\item
  $\gamma\in\cals\Leftrightarrow\phi(r(\gamma))\gamma\phi(s(\gamma))^{-1}\in
  \cals$ for every $\gamma\in\calg$ with $r(\gamma)\in A$ and $s(\gamma)\in A$. 
\end{enumerate}
Let $\aut_\cals(\calg)\vert_A\subset\emor_\cals(\calg)\vert_A$ be the subset
consisting of those $\phi$ such that $x\rightarrow r(\phi(x))$ is
a.e.\ injective. An element in $\aut_\cals(\calg)\vert_A$ is called a
\emph{local section of $\calg$}. If $A=\calg^0$, we will drop the 
subscript $A$ in the 
notation. If $\cals=\calg$, the condition ii) is void, and we will
drop the subscript $\cals$ in the notation. 
\end{definition}

\begin{definition}\label{def:composition of local sections}
Let $\phi:A\rightarrow\calg$ and $\psi:B\rightarrow\calg$ be local sections. The 
\emph{composition} $\phi\circ\psi$ is the local section defined by 
\begin{equation*}
\psi^{-1}\bigl(r(B)\cap A\bigr)\rightarrow\calg,~x\mapsto
\phi\bigl(r(\psi(x))\bigr)\circ\psi(x).
\end{equation*}
\end{definition}

\subsection{Definition of strong extensions}

The following notion is a generalization of 
(strong) normality for equivalence 
relations~\cite{feldman+sutherland+zimmer}*{Theorem~2.2 and Definition~2.14}. 

\begin{definition}\label{def:normality}
We call a subgroupoid $\cals\subset\calg$ of a discrete measured
groupoid \emph{(strongly) normal} if there is a countable 
family $\{\phi_n\}$ (in $\aut_\cals(\calg)$) in $\emor_\cals(\calg)$
such that for a.e.\ $\gamma\in\calg$ there exists exactly one $\phi_n$
in the family with $\phi_n(r(\gamma))\gamma\in\cals$. We say that
$\{\phi_n\}$ is a family of \emph{(strongly) normal choice functions}. 
\end{definition}

\begin{remark}
Note that for a translation groupoid $X \rtimes G$ and a subgroup $H \subset G$,
the groupoid $X \rtimes H \subset X \rtimes G$ is strongly normal if and only if $H$ is normal in $G$. A countable family of strongly normal choice functions is given by the automorphisms of $X$ which are induced by any complete set of representatives of the cosets $G/H$.

The previous definition makes sure that all the relevant choices can be made in a measurable way.
\end{remark}

\begin{definition}\label{def:surjection of groupoids}
Let $\theta:(\calg,\mu)\rightarrow (\calq,\nu)$ be a homomorphism 
of discrete measured
groupoids. Let $\cals=\theta^{-1}(\calq^0)$ be the kernel of $\theta$, and let
$\mu=\int_{\calq^0}\nu_yd\nu(y)$ be the measure disintegration with respect 
to $\theta$. Then 
\begin{enumerate}[i)]
\item $\theta$ is a \emph{surjection} if 
for a.e.~$x\in\calg^0$ and a.e.~$\gamma\in\calq$ with
$s(\gamma)=\theta(x)$ there exists $g\in\calg$ with $\theta(g)=\gamma$
(\emph{pointwise surjectivity}),  
\item $\theta$ is a \emph{strong surjection} if $\theta$ is a surjection 
and $(\cals\vert_{\theta^{-1}(y)}, \nu_y)$ is ergodic
for a.e.\ $y\in\calq^0$. 
\end{enumerate}
\end{definition}

The following lemma is standard and proved in the context of 
equivalence relations  
in~\cite{furman-outer}*{Lemma~2.1}. We present 
a proof of the groupoid version here so that the reader can see 
the validity of Remark~\ref{rem:parametrized version of transitivity}. 
 
\begin{lemma}\label{lem:transitivity for ergodic groupoids}
Let $\calg$ be an ergodic discrete measured groupoid. 
Let $E,F\subset\calg^0$ be Borel subsets with
the same measure.  
Then there is $\phi\in\aut(\calg)$ such that $r(\phi(E))\subset F$ and
$r\circ\phi:E\rightarrow F$ is a measure isomorphism. 
\end{lemma}

\begin{proof}
Choose a measure isomorphism 
$\calg^0\times I\rightarrow\calg$, $(x,i)\rightarrow \phi_i(x)$
as in Lemma~\ref{lem:selection theorem with injective range}. 
For any Borel subsets $A,B\subset\calg^0$ define 
\begin{equation*}
c(A,B)\defq\sup_{n\in\bbN}\mu\bigl(r(\phi_n(A))\cap B\bigr). 
\end{equation*}
By ergodicity, $c(A,B)>0$ whenever $\mu(A)>0$ and $\mu(B)>0$. 
Let $E_0=E$ and $F_0=F$. By induction we define Borel sets
$E_n,F_n\subset\calg^0$ and elements $\psi_n\in\{\phi_i;~i\in\bbN\}$
as follows: Given $E_n$ and $F_n$ let $i\in\bbN$ be the minimal number
such that
\begin{equation*}
\mu(r(\phi_i(E_n))\cap F_n)\ge c(E_n,F_n)/2, 
\end{equation*}
and let $\psi_n=\phi_i$ and 
$E_{n+1}=E_n\backslash (r\circ\phi_n)^{-1}(F_n)$, 
$F_{n+1}=F_n\backslash (r\circ\psi_n)(E_n)$. 
Let $E_\infty=\bigcap_nE_n$ and $F_\infty=\bigcap_nF_n$. 
For every $n\in\bbN$, the sets $E_n$ and $F_n$ have the 
same measure, thus $\mu(E_\infty)=\mu(F_\infty)$. Further, 
we have $\mu(E_\infty)=\mu(F_\infty)=0$ since otherwise 
$c(E_n,F_n)\ge c\defq c(E_\infty,F_\infty)>0$ for all $n$ 
which contradicts the choice of $\psi_n$ at the stage 
where $\mu(E_n\bs E_{n+1})<c/2$. 
So $(E_n')_{n\in\bbN}$ and $(F_n')_{n\in\bbN}$ defined by 
$E'_{n}=E_{n}\backslash E_{n+1}$
and $F'_{n}=F_n\backslash F_{n+1}$ are partitions of $E$ and $F$,
respectively. Now we can define $\phi:\calg^0\rightarrow\calg$ as 
follows: 
\begin{equation*}
\phi(x)=\begin{cases}
                 \id_x &\text{ if $x\not\in E$,}\\
                 \psi_n(x) &\text{ if $x\in E'_n$.}
        \end{cases}\qedhere
\end{equation*}
\end{proof}

\begin{remark}\label{rem:parametrized version of transitivity}
After fixing a measure isomorphism $\calg^0\times I\rightarrow\calg$, 
the construction of $\phi\in\aut(\calg)$ in the previous lemma 
involves no other choices. This fact is important in a situation where one 
wants to conclude that $\phi$ depends on the 
input data $\mu,E,F$ in a measurable way. 
\end{remark}

The next lemma turns out to be crucial for the 
construction of our spectral sequence in
Section \ref{sec: construction}.

\begin{lemma}\label{lem:lifting partial isomorphisms}
Let $\theta:(\calg,\mu)\rightarrow (\calq,\nu)$ be a strong surjection 
of ergodic, discrete measured groupoids. For every Borel subset 
$A\subset\calq^0$ and every $\phi\in\aut(\calq)\vert_A$ there is 
a lift $\psi\in\aut(\calg)\vert_{\theta^{-1}(A)\cap\calg^0}$, that is, 
$\theta(\psi(x))=\phi(\theta(x))$ for
a.e.\ $x\in\theta^{-1}(A)\cap\calg^0$. 
\end{lemma}

\begin{proof}
Let $\cals$ be the kernel of $\theta$.  
If the cardinality of $\calg^0$ is finite, Lemma \ref{lem:lifting partial isomorphisms} becomes very
easy and is left to the reader. 
We assume in the sequel that $\calg^0$ is infinite. Thus 
there is a Borel isomorphism $\calg^0\cong [0,1]$. Let us fix one. 
Let $\mu=\int_{\calq^0}\mu_yd\nu(y)$ be the disintegration 
of $(\calg^0, \mu)$ with respect
to $\theta$. By assumption, $\cals\vert_{\theta^{-1}(y)}$ is ergodic 
for a.e.~$y\in \calq^0$. 

First of all, we reproduce an argument in~\cite{glasner}*{Proof of
  Theorem~3.18 on p.~70}. The function 
$p:x\mapsto\mu_{\theta(x)}(\{x\})$ on $\calg^0$ is measurable and 
$\calg$-invariant. By 
ergodicity we have $p(x)=c$ for $\mu$-a.e.~$x\in\calg^0$ and  
a constant $c\ge 0$. If $c>0$, then $\theta:\calg^0\rightarrow\calq^0$
has finite fibers $\nu$-a.e., and the lemma is easily proved using the
selection theorem. We leave that to the reader and restrict ourselves
to the more complicated case that $c=0$, that is, $\mu_y$ is
continuous for $\nu$-a.e. $y\in\calq^0$. 

Consider the following Borel maps 
\begin{align*}
f:\calg^0\rightarrow \calq^0\times [0,1],&~~f(x)=(\theta(x),
\mu_{\theta(x)}([0,x])),\\ 
g:\calq^0\times [0,1]\rightarrow \calg^0,&~~(y,t)\mapsto\min\{x\in
\calg^0;\mu_y([0,x])\ge t\}. 
\end{align*}
One sees that $g\circ f=\id$. In particular, $f$ is injective. Further
we have for a Borel subset $A\subset\calq^0$  
\begin{align*}
\mu\bigl(f^{-1}(A\times [0,t])\bigr)&=
  \int_A\mu_y\bigl(g(A\times [0,t])\bigr)d\nu(y)\\
&=\int_A\mu_y\bigl(g(\{y\}\times [0,t])\bigr)d\nu(y)\\
&=\int_A\mu_y\bigl([0,g(y,t)]\bigr)d\nu(y)=t \cdot \nu(A).\\
\end{align*}
Thus $f_\ast\mu=\nu\times\lambda$, where
$\lambda$ denotes the Lebesgue measure. Hence $\Im(f)$ 
has full measure, and $f$ is a measure
isomorphism $(\calg^0,\mu)\rightarrow (\calq^0\times
[0,1],\nu\times\lambda)$. 

Since $\theta$ is pointwise surjective, 
there is a, at least set-theoretic, lift $\psi'$ of $\phi$, and 
one uses Lemma~\ref{lem:selection theorem with injective range} to
obtain a measurable such 
$\psi'\in\emor(\calg)\vert_{\theta^{-1}(A)\cap\calg^0}$. 
By the selection theorem there is a 
partition 
\[\theta^{-1}(A)\cap\calg^0=\bigcup_{k\ge 1} D_k\text{ with 
$\psi'\vert_{D_k}\in\aut(\calg)\vert_{D_k}$ for any $k\ge 1$.}\] 
Let $D'_k\defq r(\psi'(D_k))$. In the following, we neglect 
null sets. 
By the uniqueness of disintegration we obtain 
that 
\begin{equation*}\label{eq:equality of measures}
\mu_{r(\phi(y))}(D_k')=\mu_y(D_k)\text{ for every $y\in\calq^0$}.
\end{equation*}
The problem is that the sets
$\{D_k'\}_{k\ge 1}$ could have overlaps, thus causing $\psi$ to be
non-injective. The idea is to make these sets fiberwise disjoint by an
application of Lemma~\ref{lem:transitivity for ergodic groupoids} to
$(\cals,\mu_y)$ for $y\in C$. To make this precise,
consider for $y\in\calq^0$ the unique minimal sequence 
$0=m(y,0)<m(y,1)<m(y,2)<\ldots$ of numbers in $[0,1]$ with the
property 
\[\mu_{r(\phi(y))}(D_k')=\mu_{r(\phi(y))}
\bigl(g(\calq^0\times [m(y,k-1),m(y,k)])\bigr).\]
Let $\tau^y_k:D_k'\rightarrow\cals$
be the local section of $(\cals,\mu_{r(\phi(y))})$ constructed in 
Lemma~\ref{lem:transitivity for ergodic groupoids} with the
property that 
\[r\circ\tau^y_k:D_k'\rightarrow g\bigl(\calq^0\times
[m(y,k-1),m(y,k)]\bigr)\]
is a $\mu_{r(\phi(y))}$-measure isomorphism. We define
a local section $\psi$ on $\theta^{-1}(A)$ by 
\[\psi(x)=\tau^{\theta(x)}_k\circ\psi'(x)\text{ for $x\in D_k$.}\]
It follows from the explicit construction of $\tau^y_k$ (see
Remark~\ref{rem:parametrized version of transitivity}) that $\psi$ is
measurable. 
It is clear that $\psi$ is still a lift of $\phi$ since
$\tau^{\theta(x)}_k$ is (fiberwise) a local section of $\cals=\ker(\theta)$. 
One easily verifies that $\mu(r(\psi(B)))=\mu(B)$ for all 
Borel subsets $B\subset\theta^{-1}(A)$, thus 
$\psi:\theta^{-1}(A)\rightarrow\calg$ lies indeed in 
$\aut(\calg)\vert_{\theta^{-1}(A)\cap\calg^0}$. 
\end{proof}

\subsection{Strongly normal subgroupoids vs.\ strong quotients}

\begin{theorem}\label{thm:kernel is strongly normal}
Let $\theta:\calg\rightarrow\calq$ be a strong surjection of ergodic,
discrete measured groupoids. 
Then the kernel $\ker(\theta)=\theta^{-1}(\calq^0)$  
is a strongly normal subgroupoid of $\calg$. 
\end{theorem}

\begin{proof} 
Consider a measure isomorphism 
$\phi:\calq^0\times I\rightarrow\calq$ as in Theorem~\ref{lem:selection
  theorem with injective range}. According to 
Lemma~\ref{lem:lifting partial isomorphisms} we can lift each 
$\phi_i\defq\phi\vert_{\calq^0\times\{i\}}\in\aut(\calq)$ to 
$\psi_i\in\aut(\calg)$. Then 
$\{\psi_i\}_{i\in I}$ is a family of strongly normal choice
functions for $\ker(\theta)$.
\end{proof}

The following theorem is a straightforward generalization
of~\cite{feldman+sutherland+zimmer}*{Theorem~2.2}, where $\cals$ and 
$\calg$ are assumed to be equivalence relations, 
to groupoids. One does the same constructions as in~\emph{loc.~cit.} 
line by line, but for groupoids.  

\begin{theorem}[Quotient construction]\label{thm:quotient construction}
Let $\cals\subset\calg$ be a strongly normal subgroupoid of an ergodic
discrete measured groupoid. Then there
is a strong surjection $\theta:\calg\rightarrow\calq$
onto an ergodic discrete measured groupoid, called the \emph{quotient
  of $\calg$ by $\cals$}, such that 
\begin{enumerate}[i)]
\item $\ker(\theta)\defq\theta^{-1}(\calq^0)=\cals$, 
\item for a.e.~$x\in\calg^0$ and every $\gamma\in\calq$ with
  $\theta(x)=s(\gamma)$ there exists $g\in\calg$ such
  $\theta(g)=\gamma$, 
\item for any ergodic discrete measured groupoid $\calq'$
  and any homomorphism
  $\theta':\calg\rightarrow\calq'$ with $\cals\subset\ker\theta'$ 
there is a m.p.\ homomorphism
  $\kappa:\calq\rightarrow\calq'$ such that
  $\kappa\circ\theta=\theta'$. 
\end{enumerate}
\end{theorem}

\begin{definition}\label{def:strong extension}
With the setting of the previous theorem, we say that $(\calg,\mu)$ is 
a \emph{strong extension} of $(\cals,\mu)$ by $(\calq,\nu)$, and we indicate 
this, similarly to groups, by writing: 
\[1\rightarrow (\cals,\mu)\longrightarrow (\calg,\mu)
\longrightarrow (\calq,\nu)\rightarrow 1\]
\end{definition}

We record the following lemma for later reference. 

\begin{lemma}\label{lem:factorizing into lifts and kernel sections}
We retain the notation of the preceding Theorem~\ref{thm:quotient
  construction}. 
For every local section $\psi:A\rightarrow\calg$ there is a countable
Borel partition $A=\bigcup_{n\ge 1}A_n$, and local sections 
$q_n\in\aut_\cals(\calg)\vert_{A_n}$ and $s_n\in\aut(\cals)$ such that 
$\psi\vert_{A_n}=q_n\circ s_n$
\end{lemma}

\begin{proof}
Let $\{\phi_n\}_{n\ge 1}$ be a family of strongly normal choice
functions. Define 
\[A_n\defq\{a\in A;~(\phi_n\circ\psi)(a)\in\cals\}.\]
Then $s_n\defq\phi_n\circ\psi\vert_{A_n}\in\emor(\cals)\vert_{A_n}$. 
Upon further partitioning each $A_n$, we can assume that 
$s_n\in\aut(\cals)$. For $q_n=\phi_n^{-1}$ the assertion follows. 
\end{proof}

\begin{remark}\label{rem:automorphism over kernel are lifts}
Let $\cals\subset\calg$ be strongly normal with quotient
$\calg\xrightarrow{\theta}\calq$ as in Definition~\ref{def:strong extension}. 
A local section $\phi\in\aut_\cals(\calg)$ is the lift of some
$\phi'\in\aut(\calq)$ by the following argument: 
The map $\theta\circ\phi:\calg^0\rightarrow\calq$ is $\cals$-invariant. 
The restrictions of $\cals$ to fibers $\theta^{-1}(q)$ for $q\in\calq^0$ 
are ergodic since $\theta$ is a strong surjection. 
Thus, $\theta\circ\phi$ descends to a map 
$\calq^0\rightarrow\calq$ of which $\phi$ is a lift.  
\end{remark}

\begin{remark}\label{rem:lifting for Gaboriaus extension}
Consider Gaboriau's extension~\eqref{eq:Gaboriaus extension} of 
Subsection~\ref{subsec:example with groups by gaboriau} 
	\begin{equation*}
	  1\rightarrow (Z\rtimes\Lambda,\mu\times\nu)\longrightarrow 
	  (Z\rtimes\Gamma,\mu\times\nu)\longrightarrow (Y\rtimes Q,\nu)\rightarrow 1
	\end{equation*}
associated to a group extension $1\rightarrow\Lambda\rightarrow\Gamma\xrightarrow{p} Q\rightarrow 1$. 
If the action of $\Gamma$ on $Z=X\times Y$ 
is not ergodic, then this is not a strong extension 
in the sense of Definition~\ref{def:strong extension}. However, the conclusion from 
the combination of Lemma~\ref{lem:factorizing into lifts and kernel sections} and 
Remark~\ref{rem:automorphism over kernel are lifts} is still true: For every local 
section $\psi:A\rightarrow Z\rtimes\Gamma$ there is a countable Borel partition 
$A=\bigcup_{n\ge 1}A_n$, local sections $q_n:A_n\to Z\rtimes\Gamma$ 
that are 
lifts of local sections of $Y\rtimes Q$, and local sections 
$s_n:A_n\to Z\rtimes\Lambda$ such that $\psi\vert_{A_n}=q_n\circ s_n$. 

For that, note that a local section $\psi:A\rightarrow Z\rtimes\Gamma$ is essentially  
given by a map $A\rightarrow\Gamma$ that we denote by the same name. 
There are a countable Borel partition $Y=\bigcup_{n\ge 1}Y_n$ 
and elements $q_n\in Q$ for every $n\ge 1$ such that 
$p(\psi(a))=q_n$ if $a\in A_n=A\cap X\times Y_n$. Choose lifts $\gamma_n\in\Gamma$ 
for each $q_n\in Q$. Then the desired $q_n:A_n\rightarrow\Gamma$ 
and $s_n:A_n\rightarrow\Lambda$ are defined 
by $q_n(a)=\gamma_n$ and $s_n(a)=\psi(a)q_n(a)^{-1}$. 
\end{remark}

\section{Construction of the spectral sequence} \label{sec: construction}

\subsection{$\calr(\calg)$-modules}\label{subsec:modules}

Let $(\calg,\mu)$ be a discrete measured groupoid with source and
range maps $s,r\colon \calg \to \calg^0$.  Denote by $L(\calg,\mu)$
the associated von Neumann algebra and by $\calu(\calg,\mu)$ the
algebra of operators, which are affiliated with the finite von Neumann
algebra $L(\calg,\mu)$.  Recall that a local section of
$(\calg,\mu)$ is a Borel
section $\phi\colon A \to \calg$ of $s$ such that the restriction of
$r$ to $\phi(A)$ is injective.

Let $\calr(\calg,\mu)$ be the convolution algebra of complex valued
Borel functions on $\calg$, which can be written as \emph{finite} sums
of products of essentially bounded complex-valued Borel functions on
$\calg^0$ and characteristic functions $\chi_{\phi(A)}$ of 
graphs of local sections $\phi:A\rightarrow\calg$. 
Clearly, $\calr(\calg,\mu)
\subset L(\calg,\mu)$ is a sub-ring.

\begin{remark}
	We consider the case $(\calg,\mu)=(X\rtimes\Gamma,\mu)$ of a translation groupoid. The 
	\emph{crossed product ring} $L^\infty(X)\rtimes\Gamma$ is the free 
	$L^\infty(X)$-module with basis $\Gamma$. Its multiplication is determined 
	by the rule $\gamma f=f(\gamma^{-1}\_)\gamma$ 
	for $\gamma\in\Gamma, f\in L^\infty(X)$. The crossed product ring embeds 
	as a rank-dense subring into $\calr(X\rtimes\Gamma,\mu)$ via 
	$\sum f_gg\mapsto f\in L^\infty(X\rtimes\Gamma)$ 
	with $f(x,g^{-1})=f_g(x)$ (compare~\cite{thomrank}*{Proposition~4.1}). 
\end{remark}

The following lemma is a useful characterization of modules over
$\calr(\calg,\mu)$, which should be seen in analogy to group rings.

\begin{lemma}Let $M$ be an abelian group. To give a
  $\calr(\calg,\mu)$-module structure on $M$ is the same as to give a
  $L^{\infty}(\calg^0)$-module structure and compatible partial
  isomorphism as follows: For each local section $\phi$ of
  $(\calg,\mu)$, there exists an isomorphism
$$\hat\phi\colon \chi_A M \to \chi_{r(\phi(A))} M$$ that is compatible 
with composition, restriction and orthogonal sum.
\end{lemma}
\begin{proof}
  Denote by $\overline \calr(\calg,\mu)$ the universal
  $L^{\infty}(\calg^0)$-algebra, generated by symbols $\overline\phi$
  for every local section $\phi$ of $\calg$, subject to relations
  implementing compatibility with restriction, composition and
  orthogonal sum. Clearly, the natural map $\sigma\colon \overline
  \calr(\calg,\mu) \to \calr(\calg,\mu)$ is surjective and $M$ is a
  $\overline \calr(\calg,\mu)$-module. The proof is finished by
  showing that $\sigma$ is injective and hence an isomorphism.

  Let $h= \sum f_i\overline\phi_i$ be a finite sum in $\overline
  \calr(\calg,\mu)$. It follows from
  the compatibility with respect to composition and restriction that
  every element in $\overline \calr(\calg,\mu)$ can be written in this
  form. Assume that $\sum f_i \phi_i=0$. In order to show injectivity, we have to prove $h=0$.
  We may assume that the support of $f_i$ is equal to the
  domain of $\phi_i$. If we partition $\calg^0$ into sets
  $X_1,\dots,X_k$, on which
  \begin{enumerate}
  \item the domains of the local sections $\phi_i|_{X_j}$ are either
    $X_j$ or empty, and
  \item the local sections $\phi_i|_{X_j}$ are either a.e.~equal or
    a.e.~different,
  \end{enumerate}
  we can restrict our attention to one set $X_j$ at a time.  Indeed,
  since $\overline{\phi_i}$ is the orthogonal sum of its restrictions
  to the $X_j$, it suffices to show $\sum_i f_i|_{X_j} \overline
  \phi_i|_{X_j} =0$. Without loss of generality, we can now assume
  that the local sections $\phi_i|_{X_j}$ are a.e.\ different from
  each other. Clearly, $\sum _if_i|_{X_j} \phi_i|_{X_j} =0$ is only
  possible if $f_i=0$, since the $\phi_i|_{X_j}$ have disjoint support
  as characteristic functions on $\calg$. This finishes the proof.
\end{proof}

\begin{remark}
	The abelian group $L^\infty(\calg^0)$ becomes an $\calr(\calg,\mu)$-module 
	via 
	\[
		\hat{\phi}:L^\infty(A)\rightarrow L^\infty(r(\phi(A))),
		~f\mapsto \bigl(x\mapsto f((r\circ\phi)^{-1}(x))\bigr)
	\]	
	for a local section $\phi:A\rightarrow\calg$. 
\end{remark}

\subsection{Completion and localization of modules}\label{subsec:completion}
Every $L^\infty(\calg^0)$-module $M$ carries a canonical metric (\emph{rank metric}), 
which
is induced by the so-called \emph{rank}, i.e.
\[d(\xi,\eta) = \inf \{\mu(A^c) \mid A \subset \calg^0 \mbox{ Borel}, 
\chi_A\xi = \chi_A\eta\}, \quad \forall \xi,\eta \in M.\]

It was shown in \cite{thomrank}*{Lemma 4.4} that the completion of a
$\calr(\calg,\mu)$-module with respect to the underlying
$L^\infty(\calg^0)$-module carries a natural $\calr(\calg,\mu)$-module
structure and that the associated completion functor is exact~\cite{thomrank}*{Lemma 2.6}.  Moreover, the category of complete
$\calr(\calg,\mu)$-modules was shown to be abelian with enough
projective objects~\cite{thomrank}*{Theorem 2.7}.  

We denote by $\compmod{\calr(\calg,\mu)}$ the full subcategory of 
$\modules{\calr(\calg,\mu)}$, formed by complete $\calr(\calg,\mu)$-modules and
denote by
$$c\colon \modules{\calr(\calg,\mu)} \to \compmod{\calr(\calg,\mu)}$$
the completion functor.

Projectives in $\compmod{\calr(\calg,\mu)}$ 
are obtained by completing free modules.  The following
lemma shows that there are also enough injective objects.

\begin{lemma}
  The abelian category $\compmod{\calr(\calg,\mu)}$ has enough
  injective objects.
\end{lemma}
\begin{proof}
  Clearly, the abelian category $\modules{{\calr(\calg,\mu)}}$ has
  enough injective objects.  Let $Z$ be an injective object in the
  abelian category $\modules{\calr(\calg,\mu)}$. It is sufficient to
  construct an complete injective object $Z'$, which contains $Z$ as a
  submodule.  Let $\alpha$ be an ordinal with uncountable cofinality,
  i.e.\ there exists no countable cofinal subset.  For every $\beta
  \leq \alpha$, we define $Z_{\beta}$ by the following transfinite
  inductive procedure. We set $Z_0 = Z$ and $Z_{\beta+1}$ to be some
  injective ${\calr(\calg,\mu)}$-module, containing the completion of
  $Z_{\beta}$ as a sub-module. If $\beta$ is a limit ordinal, we set
  $Z_{\beta} = \cup_{\beta' < \beta} Z_{\beta'}$. Every Cauchy
  sequence in $Z_{\alpha}$ is contained in some $Z_{\beta}$ for
  $\beta< \alpha$ and thus its limit exists in $Z_{\beta+1} \subset
  Z_{\alpha}$. Hence $Z_{\alpha}$ is complete. Moreover, choosing
  $\alpha$ big enough, $Z_{\alpha}$ is injective as a
  ${\calr(\calg,\mu)}$-module, since injectivity can be tested on the
  \emph{set} of sub-modules of the trivial module. We can define $Z' =
  Z_{\alpha}$.
\end{proof}

$\calr(\calg,\mu)$-modules, which are zero-dimensional as
$L^{\infty}(\calg^0)$-modules, complete to the zero module. Indeed,
this is a reformulation of the local criterion of zero-dimensionality,
that can be found in \cite{sauer(2005)}*{Theorem 2.4}. Hence, there
exists a natural exact functor of abelian categories
\begin{equation}\label{eq:completion functor}
  c\colon \locmod{\calr(\calg,\mu)} \to \compmod{\calr(\calg,\mu)},
\end{equation}
which is an equivalence of abelian categories. Indeed, the inverse is
the restriction of the natural quotient functor.  In the sequel we
concentrate on complete $\calr(\calg,\mu)$-modules. The completion of
$L^{\infty}(\calg^0)$ identifies naturally with the algebra of
measurable functions on $\calg^0$, which we denote by
$\calm(\calg^0)$.  Note that $\calu(\calg,\mu)$ is complete, whereas
$L(\calg,\mu)$ is not necessarily.

\subsection{Some derived functors}

In this subsection, we introduce some derived functors, which are
appropriate to make our approach to a spectral sequence work. They are
defined on the localized category; and thus we are heavily using its
nice homological properties which were established in preceding
section and \cite{thomrank}.

\begin{definition}
  Let $M$ be a complete left $\calr(\calg,\mu)$-module. We define:
$$H_*(\calg,M) = 
\underline\Tor_*^{\calr(\calg,\mu)}\left(\calm(\calg^o),M\right),\quad \mbox{and}
\quad H^*(\calg,M) =
\underline\Ext^*_{\calr(\calg,\mu)}(\calm(\calg^0),M),$$ to be the
\emph{complete} $\calg$-homology and $\calg$-cohomology of the
$\calr(\calg,\mu)$-module $M$.
\end{definition}

Here, $\underline\Tor^{\calr(\calg,\mu)}_*(\calm(\calg^0),?)$ and
$\underline\Ext^*_{\calr(\calg,\mu)}(\calm(\calg^0),?)$ denote the
derived functors of the functors
$$M \mapsto \calm(\calg^o) \otimes_{\calr(\calg,\mu)} M, \quad \mbox{and} \quad 
M \mapsto \hom_{\calr(\calg,\mu)}(\calm(\calg^0),M),$$ from the
abelian category $\compmod{\calr(\calg,\mu)}$ to abelian
groups. Following the arguments in \cite{weibel}*{Chapter 2.7}, we see
that the bi-functors $\underline\Tor^{\calr(\calg,\mu)}_*(?,?)$ and
$\underline\Ext_{\calr(\calg,\mu)}^*(?,?)$ are balanced, as are their
classical counterparts.
\subsection{Spectral sequence}
The following theorem enables us to construct a 
spectral sequence for the cohomology which we just defined.

\begin{theorem}\label{thm:spectral sequence in degree zero}
  Let $$1 \to (\cals,\mu') \to (\calg,\mu) \to (\calq,\nu) \to 1$$ be
  a strong extension of discrete measured groupoids, i.e.\
  $(\calg,\mu) \to (\calq,\nu)$ is a strong surjection with kernel
  $(\cals,\mu')$.  Let $M$ be a complete left-module over
  $\calr(\calg,\mu)$. The abelian group
$$\hom_{\calr(\cals,\mu')} (\calm(\calg^0),M)$$
is naturally a complete left module over $\calr(\calq,\nu)$, and there
is a natural isomorphism of abelian groups
\[
\hom_{\calr(\calq,\nu)}
\left(\calm(\calq^0),\hom_{\calr(\cals,\mu')}(\calm(\calg^0),M)\right)
= \hom_{\calr(\calg,\mu)}(\calm(\calg^0),M).\]
\end{theorem}
\begin{proof}
  Since $\calm(\calq^0) \subset \calm(\calg^0)$, there is a natural
  $L^{\infty}(\calq^0)$-module structure on $$\hom_{\calr(\cals,\mu')}
  (\calm(\calg^0),M).$$

  It remains to provide partial isomorphisms as described above.  Let
  $\phi$ be a local section of $\calq$. Since $(\calg,\mu) \to
  (\calq,\nu)$ is a strong surjection, there exists a section $\phi'$
  of $\calg$ which lifts $\phi$ by Lemma \ref{lem:lifting partial
    isomorphisms}.  For $f \in \hom_{\calr(\cals,\mu')}
  (\calm(\calg^0),M)$, we define 
$$(\phi \rhd f)(g) = \phi' f(\phi'^{-1}g),\quad \forall g \in \calm(\calg^0).$$
This is well-defined and, together with the aforementioned
$\calm(\calq^0)$-module structure, it defines a left-module structure
of $\calr(\calq,\nu)$. Moreover, $\hom_{\calr(\cals,\mu')}
(\calm(\calg^0),M)$ is easily seen to be complete. We now aim to show
that
\[
\hom_{\calr(\calq,\nu)}
\left(\calm(\calq^0),\hom_{\calr(\cals,\mu')}(\calm(\calg^0),M)\right)
= \hom_{\calr(\calg,\mu)}(\calm(\calg^0),M).\]

An element in $f \in \hom_{\calr(\calg,\mu)}(\calm(\calg^0),M)$ is
described by the image of $1 \in \calm(\calg^0)$. For this, note that all $L^{\infty}(\calq^0)$-module
maps are automatically contractive and hence continuous in the rank metric, 
see \cite{thomrank}*{Lemma 2.3}.

Let $\phi$ be a local section of $(\calg,\mu)$. We set $\hat \phi =
\chi_{\phi(A)} \in \calr(\calg,\mu)$.  The image of $1$ under $f$
satisfies
$$\hat \phi f(1) = f(\chi_{r(\phi(A))}) = \chi_{r(\phi(A))} f(1).$$
Conversely, every such element gives rise to a module homomorphism.

An element in $\hom_{\calr(\cals,\mu')}(\calm(\calg^0),M)$ is
described as an element in $M$, satisfying the above relation for
local sections of $(\cals,\mu')$. Having this description, it is
obvious that
$$\hom_{\calr(\calq,\nu)}
\left(\calm(\calq^0),\hom_{\calr(\cals,\mu')}(\calm(\calg^0),M)\right)$$
identifies with $\hom_{\calr(\calg,\mu)}(\calm(\calg^0),M)$ since we
can write each local section of $(\calg,\mu)$ as a countable orthogonal sum of
local sections, which are products of lifts of a local section of
$(\calq,\nu)$ and a local section of $(\cals,\mu')$ by
Lemma~\ref{lem:factorizing into lifts and kernel sections} and
Remark~\ref{rem:automorphism over kernel are lifts}. This finishes the
proof, since $M$ is complete and the orthogonal decomposition gives
rise to a rank convergent sum $M$.
\end{proof}

The decomposition of functors which was established in the preceding
theorem yields a Grothendieck spectral sequence by Theorem
\ref{grothendieck}:
\begin{theorem}\label{thm:spectral sequence}
 Let 
\begin{equation*}
1 \to (\cals,\mu') \to (\calg,\mu) \to (\calq,\nu) \to 1
\end{equation*}
be a strong extension of discrete measured groupoids.   
Let $M$ be a complete left $\calr(\calg,\mu)$-module. Then there is
  a first quadrant spectral sequence with $E_2$-term
  \begin{eqnarray}
    E_2^{pq} &=& 
        \underline\Ext^p_{\calr(\calq,\nu)}\left(\calm(\calq^0),
        \underline\Ext^q _{\calr(\cals,\mu')}
      \left(\calm(\calg^0),M \right)\right) \\ \nonumber
    &=& H^p(\calq, H^q(\cals,M))
  \end{eqnarray}
  converging to
  $\underline\Ext^{p+q}_{\calr(\calg,\mu)}\left(\calm(\calg^0),M\right)
  = H^{p+q}(\calg,M)$. For $M=\calu(\calg,\mu)$, which is a
  $\calr(\calg,\mu)$-$\calu(\calg,\mu)$-bimodule, we obtain a spectral
  sequence of $\calu(\calg,\mu)$-modules.
\end{theorem}
\begin{proof} The only claim that remains to be verified is that $\hom_{\calr(\cals,\mu')}(\calm(\calr^0),?)$
sends injective objects to injective objects. However, this 
follows from Lemma \ref{inj-criterion} if we can provide an exact left adjoint 
functor. Clearly, a similar argument as above shows that
$$M \mapsto c\left(\calm(\calr^0) \otimes_{\calm(\calq^0)} M\right)$$
is left adjoint to $\hom_{\calr(\cals,\mu')}(\calm(\calr^0),?)$. Moreover, it is the composition of a flat ring extension
(note that $\calm(\calq^0)$ is von Neumann regular) and the exact completion functor, 
and thus exact. This finishes the proof.
\end{proof}

\begin{remark}\label{rem:E1 term}
  For later reference we note that the $E_1$ is of the spectral
  sequence is given by
$$ E_1^{p,q}=\hom_{\calr(Q,\nu)}\left( P_p,\underline\Ext^q _{\calr(\cals,\mu')}
  \left(\calm(\calg^0),M \right)\right) $$ for a projective resolution
$P_\ast$ of $\calm(\calq^0)$ in $\compmod{\calr(\calq,\nu)}$.
\end{remark}

\begin{remark}\label{rem:gaboriau also works for spectral sequence}
	The conclusion of Theorem~\ref{thm:spectral sequence} also holds true for 
	Gaboriau's extension~\eqref{eq:Gaboriaus extension} even if the group $\Gamma$ does 
	not act ergodically and so Gaboriau's extension is not a strong extension in the 
	sense of Definition~\ref{def:strong extension}: 
	For that, note that the conclusion of 
	Lemma~\ref{lem:lifting partial isomorphisms}, which is used in 
	Theorem~\ref{thm:spectral sequence in degree zero}, is very easy to see for 
	Gaboriau's extension (compare the argument in 
	Remark~\ref{rem:lifting for Gaboriaus extension}). 
	Furthermore, one replace the application 
	of Lemma~\ref{lem:factorizing into lifts and kernel sections} and 
	Remark~\ref{rem:automorphism over kernel are lifts} in the proof of 
	Theorem~\ref{thm:spectral sequence in degree zero} by 
	Remark~\ref{rem:lifting for Gaboriaus extension}. All other steps in the 
	proofs of Theorems~\ref{thm:spectral sequence in degree zero} 
	and~\ref{thm:spectral sequence} stay the same in the case of Gaboriau's extension. 
\end{remark}

\subsection{Identification of $L^2$-Betti numbers}
\label{subsec:identification}

In order to apply the spectral sequence from Section \ref{sec:
  construction} to the computation of $L^2$-Betti numbers of groups,
we need to study the special case of a m.p.\ group action more
closely.  Let $\Gamma$ be a discrete group, and let $(X,\mu)$ be a
standard Borel probability space on which $\Gamma$ acts by m.p.\ Borel
isomorphisms. Note that we do \emph{not} impose any conditions like
ergodicity or freeness of the action.

In many situations, invariants of a m.p.\ action of a discrete group
are actually invariants of the group itself. In the sequel we want to
identify the dimensions of the homological and cohomological
invariants for $(X \rtimes \Gamma,\mu)$, which we just introduced,
with ordinary $L^2$-Betti numbers of $\Gamma$.

First of all, for any discrete measured groupoid $(\calg,\mu)$, there
is a natural transformation
\begin{equation} \label{eq: eq10}
  \Tor^{\calr(\calg,\mu)}_*(?,L(\calg,\mu)) \to \underline
  \Tor^{\calr(\calg,\mu)}_*(?,\calu(\calg,\mu)),
\end{equation}
consisting of dimension isomorphisms. Indeed, it is induced by the natural map
$$? \otimes_{\calr(\calg,\mu)} 
L(\calg,\mu) \to c(?) \otimes_{\calr(\calg,\mu)} \calu(\calg,\mu).$$
For a proof that this is a dimension isomorphism, we refer to the proof of \cite{thomrank}*{Proposition 4.7}.
If follows from \cite{thomrank}*{Lemma 1.1} that the induced
map on derived functors are dimension isomorphisms as well.

Secondly, there is an isomorphism of right $\calu(\calg,\mu)$-modules
$$\hom_{\calu(\calg,\mu)} 
\left(\underline \Tor^{\calr(\calg,\mu)}_*(?,\calu(\calg,\mu)),\calu(\calg,\mu) \right)
\cong \underline \Ext_{\calr(\calg,\mu)}^*(?,\calu(\calg,\mu)),$$
since $\calu(\calg,\mu)$ is self-injective. Moreover, since dualizing
is dimension preserving, we get:
\begin{equation} \label{eq10b} \dim \underline
  \Tor^{\calr(\calg,\mu)}_*(?,\calu(\calg,\mu)) = \dim \underline
  \Ext_{\calr(\calg,\mu)}^*(?,\calu(\calg,\mu)).
\end{equation}

In~\cite{sauer(2005)} the first author observed that due to dimension
flatness of the ring extensions
$$\bbC \Gamma \subset \calr(X \rtimes \Gamma),
\quad \mbox{and} \quad L(\Gamma) \subset L(X \rtimes \Gamma),$$
there are natural dimension isomorphisms as follows:
\begin{equation} \label{eq11}
\Tor^{\calr(\calg,\mu)}_*(L^{\infty}(\calg^0),L(\calg,\mu))
  \cong H_*(\Gamma,L(X \rtimes \Gamma)) \cong H_*(\Gamma,L( \Gamma))
  \otimes_{L\Gamma} L(X \rtimes \Gamma).
\end{equation}

Thus, combining the dimension isomorphisms~\eqref{eq: eq10} for
$L^{\infty}(\calg^0)$, and~\eqref{eq10b}, and~\eqref{eq11}, we get:
\begin{align}\label{eq:betti identification}
  \betti_*(\Gamma) &= 
  \dim_{L(X \rtimes \Gamma)} H_*(X \rtimes \Gamma, L(X \rtimes \Gamma))\\
  &=
  \dim_{L(X \rtimes \Gamma)} H^*(X \rtimes \Gamma, \calu(X \rtimes \Gamma))\nonumber
\end{align}

\begin{remark} \label{vanish} Note that $H^n(\calg,\calu(\calg,\mu))$
  is the $\calu(\calg,\mu)$-dual of $H_n(\calg,\calu(\calg,\mu))$ and
  thus has the pleasant feature that it vanishes as soon as its
  dimension is zero. This follows from the results in \cite{thoml2}.
\end{remark}

\section{Examples of groupoid extensions}\label{sec:examples}

\subsection{Group extensions} \label{subsec:example with groups} Let
$\Gamma$ be a group and let $\Lambda \subset \Gamma$ be a
normal subgroup. Let $(X,\mu)$ be a standard Borel space with a probability
measure $\mu$, on which $\Gamma$ acts by m.p.\ Borel isomorphisms. For
example, one can take any measure on $\{0,1\}$ and consider the
infinite product $\{0,1\}^{\Gamma}$ on which $\Gamma$ acts by shifts.

Clearly, the translation groupoid $(X \rtimes \Gamma,\mu)$ is a
discrete measured groupoid and $(X \rtimes \Lambda,\mu)$ is a strongly
normal subgroupoid. If $(X \rtimes \Gamma,\mu)$ is ergodic, then 
let $\calq^X_{\Lambda \subset \Gamma}$ denote the
quotient groupoid, which exists by Theorem \ref{thm:quotient construction}. 

\subsection{Gaboriau's extension}
\label{subsec:example with groups by gaboriau}
Let $1\rightarrow\Lambda\rightarrow\Gamma\xrightarrow{p} Q\rightarrow
1$ be a short exact sequence of groups. We describe 
a special case of~\ref{subsec:example with groups} 
that Gaboriau used to prove
vanishing results for $L^2$-Betti numbers of
groups~\cite{gaboriau(2002b)}.  

\smallskip

Let $(X,\mu)$ be a $\Gamma$-probability space and $(Y,\nu)$ be a $Q$-probability 
space. Let $Z=X\times Y$ be the product of probability spaces. 
The group $\Gamma$ acts on $Y$ via $p$ and on $Z$ by
the diagonal action. Then $(Z\rtimes\Lambda,\mu\times\nu)$ 
is a strongly normal subgroupoid of $(Z\rtimes\Gamma,\mu\times\nu)$. 
We refer to 
\begin{equation}\label{eq:Gaboriaus extension}
  1\rightarrow (Z\rtimes\Lambda,\mu\times\nu)\longrightarrow 
  (Z\rtimes\Gamma,\mu\times\nu)\longrightarrow (Y\rtimes Q,\nu)\rightarrow 1
\end{equation}
as \emph{Gaboriau's extension}; it 
is a strong extension provided $\Gamma$ acts ergodically on $Z$. 
For every ergodic $Q$-probability space $(Y,\nu)$ 
we can find an ergodic $\Gamma$-probability 
space $(X,\mu)$ such that $\Gamma$ acts ergodically on $(Z,\mu\times\nu)$: 
Take, for example, the infinite product $\{0,1\}^\Gamma$
with the equidistribution on $\{0,1\}$ and $\Gamma$ acting by the
shift (\emph{Bernoulli action}). Since this action is mixing, 
the diagonal $\Gamma$-action on $Z$ is still ergodic. 
\subsection{The principal extension}
\label{subsec:principal extension}
Let $(\calg,\mu)$ be an ergodic discrete measured groupoid. We denote by
$$\calg_{\rm stab}\defq\{\gamma \in \calg \mid r(\gamma)=s(\gamma)\}$$ 
the \emph{stabilizer groupoid}. We call $\calg_x\defq\{\gamma\in\calg;
r(\gamma)=s(\gamma)=x\}$ the \emph{isotropy group}, or the
\emph{stabilizer}, of $x\in\calg^0$.  $(\calg_{\rm stab},\mu)$ is a
strongly normal subgroupoid of $(\calg,\mu)$. Indeed, the Borel map
$(r\times s)\colon (r \times s)^{-1}(\Delta(\calg^0)) \to \calg^0$ has
countable fibers. Hence, we can use the selection theorem to find
choice functions.  We denote by $(\calg_{\rm rel},\mu)$ the quotient
groupoid and call it the \emph{groupoid of the associated equivalence
  relation}. Obviously, $(\calg_{\rm rel},\mu)$ has the same unit
space $\calg^0$.  The strong extension
$$0 \to (\calg_{\rm stab},\mu) \to (\calg,\mu) \to (\calg_{\rm rel},\mu) \to 0$$
is called the \emph{principal extension}, associated to the discrete
measured groupoid $(\calg,\mu)$.

The discrete measured groupoid $\calg_{\rm stab}$ can be best viewed
as a \emph{direct integral} or \emph{Borel field} of discrete groups. All functorial
constructions can be carried out fiberwise. The following lemma is 
therefore no surprise and we omit its proof, which needs a little technical detour.

\begin{lemma}\label{lem: betti number of groupoid with only isotropy}
With the notation of the previous example, for every $p\ge 0$ we have 
$$ \betti_p(\calg_{\rm
  stab})=\int_{\calg^0}\betti_p(\calg_x)d\mu(x). $$
\end{lemma}

\section{Measurable cohomological dimension}\label{singer}
We are now coming to a more conceptional study of the cohomological properties of a discrete measured groupoid. The main application of the results in this section is the proof of the Hopf-Singer Conjecture for poly-surface groups.

\subsection{Definition of the measurable cohomological dimension}
Let $(\calg,\mu)$ be a discrete measured groupoid and $\compmod{\calr(\calg,\mu)}$ the category of complete $\calr(\calg,\mu)$-modules. We define the \emph{measurable cohomological dimension} ${\rm mcd}_{\bbC}(\calg,\mu)$ to be the length of the shortest resolution of $\calm(\calg^0)$ by projective objects in $\compmod{\calr(\calg,\mu)}$.
The following result is classical~\cite{weibel}*{Lemma~4.1.6 on p.~93}:

\begin{theorem} Let $(\calg,\mu)$ be a discrete measured groupoid. The following statements are equivalent:
\begin{enumerate}[i)]
\item ${\rm mcd}_{\bbC}(\calg,\mu) \leq n$.
\item For all $M \in \compmod{\calr(\calg,\mu)}$ and $m >n$, one has
$$\underline\Ext^m _{\calr(\calg,\mu)}   \left(\calm(\calg^0),M \right) =0.$$
\end{enumerate}
\end{theorem}
The spectral sequence of Theorem~\ref{thm:spectral sequence} 
immediately yields a product estimate as follows:
\begin{corollary}[of Theorem~\ref{thm:spectral sequence} and Remark~\ref{rem:gaboriau also works for spectral sequence}]\label{cor:sub-additivity}
 Let $$1 \to (\cals,\mu') \to (\calg,\mu) \to (\calq,\nu) \to 1$$ be
  a strong extension of discrete measured groupoids or 
a Gaboriau's extension~\eqref{eq:Gaboriaus extension}. Then,
  $${\rm mcd}_{\bbC}(\calg,\mu) \leq {\rm mcd}_{\bbC}(\cals,\mu') + {\rm mcd}_{\bbC}(\calq,\nu).$$ 
\end{corollary}

Turning to the interesting case of translation groupoids, we make the following definition:
\begin{definition}
Let $\Gamma$ be a discrete group. We set:
$${\rm mcd}_{\bbC}(\Gamma) = \min\{{\rm mcd}_{\bbC}(X \rtimes \Gamma,\mu) \mid \Gamma \curvearrowright (X,\mu) \},$$
where we take the minimum over all m.p. actions of $\Gamma$ 
on a standard probability space $(X,\mu)$. We call 
${\rm mcd}_{\bbC}(\Gamma)$ 
the \emph{measurable cohomological dimension} of $\Gamma$.
\end{definition}

We need the following lemma.
\begin{lemma} \label{techlem}
Let $\Gamma$ be a discrete group and let $\Gamma \curvearrowright (X,\mu)$ and $\Gamma \curvearrowright (Y,\nu)$ measure preserving actions as above. Let $f\colon Y \to X$ be a measure preserving $\Gamma$-equivariant map. Then,
$${\rm mcd}_{\bbC}(Y \rtimes \Gamma,\nu) \leq {\rm mcd}_{\bbC}(X \rtimes \Gamma,\mu).$$
\end{lemma}
\begin{proof}
Clearly, there is a $\Gamma$-equivariant ring homomorphism $f^*\colon \calm(X) \to \calm(Y)$.
Let $P_* \to \calm(X)$ be a projective resolution in the category $\compmod{\calr(X \rtimes \Gamma,\mu)}$. Since $\calm(X)$ is von Neumann regular, $\calm(Y)$ is a flat
$\calm(X)$-module. Hence,
$$\left( \calm(Y) \rtimes\Gamma \right) \otimes_{\calm(X) \rtimes\Gamma} P_* = 
\calm(Y) \otimes_{\calm(X)} P_* \to \calm(Y) \otimes_{\calm(X)} \calm(X) = \calm(Y)$$
is a resolution. Setting $Q_*=\left( \calm(Y) \rtimes\Gamma \right) \otimes_{\calm(X) \rtimes\Gamma} P_*$ and using exactness of completion (with respect to the $\calm(Y)$-module structure), we obtain a resolution $Q'_* \to \calm(Y)$ by complete modules. 
The modules $Q'_*$ are projective in 
$\compmod{\calr(Y \rtimes \Gamma,\nu)}$~\cite{thomrank}*{Theorem~2.7~(3)}. 
\end{proof}

\begin{proposition}\label{pro:subadditivity for groups}
Let $1 \to \Lambda \to \Gamma \to Q \to 1$ be an extensions of groups. Then
$${\rm mcd}_{\bbC}(\Gamma) \leq {\rm mcd}_{\bbC}(\Lambda) + {\rm mcd}_{\bbC}(Q).$$
\end{proposition}

\begin{proof}
Let $\Lambda \curvearrowright (X',\mu')$ and $Q \curvearrowright (Y,\nu)$ be measurable 
actions as above that achieve the measurable cohomological dimension.
Let $X$ be the coinduction
  of $X'$, i.e. the $\Gamma$-space $X=\map_\Lambda(\Gamma,X')$, 
on which $\gamma\in\Gamma$ acts
  from the left by composition with the $\Gamma$-map $r_{\gamma^{-1}} \colon \Gamma
  \rightarrow\Gamma,~ \gamma_0 \mapsto \gamma_0\gamma^{-1}$. 
By choosing a set theoretic section
  $s \colon \Gamma/\Lambda\rightarrow\Gamma$ of the projection with $s(1)=1$, we obtain a
  bijection
  $X \xrightarrow{\cong} \prod_{\Gamma/\Lambda} X'$. 
  We endow $X$ with
  the structure of a standard probability space $(X,\mu)$ by 
  pulling back the product measure.  This
  structure does not depend on the choice of $s$; the measure $\mu$ is
  $\Gamma$-invariant~\cite{gaboriau-examples}*{3.4}. 
We have  
${\rm mcd}_{\bbC}(X \rtimes \Lambda, \mu) \leq {\rm mcd}_{\bbC}(X' \rtimes \Lambda,\mu')$ by the previous lemma. Setting $Z = X \times Y$ as in Gaboriau's extension (see Subsection \ref{subsec:example with groups by gaboriau}) yields a translation groupoid $(Z \rtimes \Gamma, \mu \times \nu)$ and an extension of discrete measured groupoids as follows:
\begin{equation*}
  1\rightarrow (Z\rtimes\Lambda,\mu\times\nu)\longrightarrow 
  (Z\rtimes\Gamma,\mu\times\nu)\longrightarrow (Y\rtimes Q,\nu)\rightarrow 1
\end{equation*}
Again by the previous lemma, 
${\rm mcd}_{\bbC}(Z \rtimes \Lambda, \mu\times\nu) \le{\rm mcd}_{\bbC}(X \rtimes \Lambda, \mu)$. By Corollary~\ref{cor:sub-additivity} we obtain that 
$${\rm mcd}_{\bbC}(\Gamma) \leq {\rm mcd}_{\bbC}(X' \rtimes \Lambda,\mu') + {\rm mcd}_{\bbC}(Y \rtimes Q,\nu).$$ Hence we completed the proof.
\end{proof}

The precise relationship between the measurable cohomological dimension to Gaboriau's ergodic dimension~\cite{gaboriau(2002b)} is not clear at present. Certainly, the measurable cohomological dimension is smaller or equal to the ergodic dimension; however, the reverse inequality seems to be more difficult to establish.

\subsection{The Singer condition}\label{sub:singer condition}

Striving for a conceptional explanation of the Hopf-Singer conjecture in terms of ergodic theory, we make the following definition. Recall that 
a \emph{Poincar\'e duality group of dimension $n$} 
is a group $\Gamma$ of cohomological dimension $n$ such that 
\[H^p(\Gamma,\bbZ \Gamma)\cong \begin{cases} \bbZ^\epsilon &\text{ if $p=n$,}\\ 
 											  0 & \text{ if $p\ne n$.} 
								\end{cases}\]
Here $\bbZ^\epsilon=\bbZ$ as an abelian group. 
The group 
$H^p(\Gamma,\bbZ \Gamma)$ carries a natural right $\Gamma$-module structure, and 
$\bbZ^\epsilon$ denotes the right $\Gamma$-module structure on $\bbZ$ that makes 
the above isomorphism $\Gamma$-equivariant. If $\bbZ^\epsilon=\bbZ$ is the trivial 
$\Gamma$-module, then $\Gamma$ is called \emph{orientable}. 

\begin{definition} \label{defsinger}
Let $\Gamma$ be a Poincar\'e duality group of dimension $2n$. 
We say that $\Gamma$ satisfies \emph{Singer's condition} if
$$ {\rm mcd}_{\bbC}(\Gamma) \leq n.$$
\end{definition}

\begin{remark}
The definition could be phrased in a straightforward way for proper and cocompact $\Gamma$-CW-complexes, satisfying equivariant Poincar\'e duality. We leave this to the reader.
\end{remark}

Note that Singer's condition and~\eqref{eq:betti identification} 
imply that all $L^2$-Betti numbers above the middle dimension vanish. 
Hence, by Poincar\'e duality, the only possible 
non-trivial $L^2$-Betti number is in the middle dimension. Note that this is also true in the non-orientable case since every 
non-orientable Poincare duality group has an orientable subgroup of index $2$. 
It is thus of some interest to understand the class of groups which satisfy Singer's condition. Fundamental groups of closed surfaces with genus $\ge 2$ 
are measure equivalent to the free 
group of rank~$2$~\cite{gaboriau-examples}. The free abelian 
groups $\bbZ^n$ are measure equivalent to $\bbZ$~\cite{gaboriau-examples}. 
Because of Lemma~\ref{newlemma} we thus obtain the first interesting examples 
of groups satisfying Singer's condition: 

\begin{theorem}
Let $S_g$ be the closed orientable surface of genus $g\geq 1$. Then 
$\pi_1(S_g)$ satisfies Singer's condition.
\end{theorem}

Since Poincare duality groups are closed under extensions~\cite{bieri}*{Satz~2.5}, 
Proposition~\ref{pro:subadditivity for groups} yields: 

\begin{theorem}
Let $\Lambda$ be a Poincar\'e duality group of dimension $2n$ and $Q$ 
be a Poincar\'e duality group of dimension $2m$. Let 
$1 \to \Lambda \to \Gamma \to Q \to 1$
be an extension of groups. Then $\Gamma$ is a Poincar\'e duality group 
of dimension $2(n+m)$, 
and it satisfies Singer's condition if $\Lambda$ and $Q$ satisfy Singer's condition.
\end{theorem}

It is obvious that this provides a powerful tool for the construction of groups which satisfy Singer's condition. Indeed, free products of amenable groups and surface groups (in fact all groups measure equivalent to free groups) can be used as basic building blocks. A class of groups which has been studied to some extend and fits into our framework is given by the following definition:

\begin{definition} \label{defpolysurf}
A group $\Gamma$ is said to be a \emph{poly-surface group}, if there exists
a series of normal subgroups
$$\{e\} = \Gamma_0 \lhd \Gamma_1 \lhd \cdots \lhd \Gamma_n = \Gamma,$$
such that $\Gamma_k/\Gamma_{k-1}$ is a surface group, i.e.~the fundamental group 
of a closed, orientable surface of genus $\ge 2$, for all $1 \leq k \leq n$.
\end{definition}

It was proved in \cite{polysurface} that every poly-surface group admits a subgroup of finite index, which is the fundamental group of a closed, orientable, 
aspherical manifold. Hence, the following corollary is of importance.

\begin{corollary} \label{polysinger}
Poly-surface groups satisfy Singer's condition. In particular, the Hopf-Singer conjecture holds for closed aspherical manifolds with poly-surface fundamental group.
\end{corollary}

\section{Proofs of applications}\label{sec:applications}

\subsection{Proof of Theorem~\ref{thm:vanishing with amenable quotient}}
We give the proof only in the case where $\calq=\calr_{\hyp}$ is an infinite amenable  equivalence relation. (The more general case of infinite amenable measured groupoids 
is more tedious but analogous.)
There exists an increasing chain of finite sub-relations
$\calr^n_{\hyp} \subset \calr_{\hyp}$, which are isomorphic to the groupoid $(X \rtimes S_n)_{\rm rel}$, where
$S_n$ permutes a partition of a continuous Borel probability space $(X,\mu)$ into $n$ sets of equal measure. Moreover, we find these sub-relations such
that
$$\calr_{\hyp} = \bigcup_{n=1}^{\infty} \calr^n_{\hyp}.$$
In the extension
$$ 1 \to \cals \to \calr \stackrel{\phi}{\to} \calr_{\hyp} \to 1$$
we can take inverse images $\calr^n = \phi^{-1}(\calr^n_{\hyp})$ and note that
$$n \cdot \betti_p(\calr^n) = \betti_p(\cals).$$
Indeed, this follows from standard homological arguments  
since $\calr(\calr^n,\mu)$ is just an $n\times n$-matrix algebra over $\calr(\cals,\mu)$.
Since $\calr = \cup_{n=1}^{\infty} \calr^n$, we conclude that
$$\betti_p(\calr) \leq \liminf_{n \to \infty} \betti_p(\calr^n) = 
\liminf_{n \to \infty} n^{-1} \cdot \betti_p(\cals).$$
Provided the first inequality holds, this would finish the proof
 since $\betti_p(\cals)<\infty$ implies now
$\betti_p(\calr)=0$. The first inequality follows from the following three facts:

\begin{enumerate}[(i)]
\item \label{exact} $\calr(\calr^n,\mu) \subset \calr(\calr,\mu)$ 
is a dimension flat ring extension since $\cup_{n} \calr(\calr^n,\mu)$ is flat over $\calr(\calr^n,\mu)$ and rank-dense
in $\calr(\calr,\mu)$,
\item \label{eqcolim} \begin{equation*} 
L^{\infty}(\calr^0) \otimes_{\calr(\calr,\mu)} ? = 
c\left(\colim L^{\infty}(\calr^0) \otimes_{\calr(\calr^n,\mu)} ? \right) 
\end{equation*} 
for the same reason that 
$\cup_{n} \calr(\calr^n,\mu)$ is rank-dense in $\calr(\calr,\mu)$, and
\item \label{colimdim}
$$\dim \colim M_n \leq \liminf_{n \to \infty} \dim M_n$$
by \cite{lueck(2002)}*{Theorem~6.13~(2) on p.~243}. 
\end{enumerate}

By \eqref{exact}, \eqref{eqcolim}, and exactness of completion and directed co-limits,
\begin{eqnarray*}
H_k(\calr,\calu(\calr,\mu)) &=& \colim H_k(\calr^n,\calu(\calr^n,\mu)) \\
&=& \colim \left(H_k(\calr^n,\calu(\calr^n,\mu)) \otimes_{\calu(\calr^n,\mu)}
\calu(\calr,\mu) \right).
\end{eqnarray*}
This implies the claim since, using \eqref{colimdim}, we get:
\begin{eqnarray*}
\betti_k(\calr^n) &=&
\dim_{\calu(\calr^n,\mu)} H_k(\calr^n,\calu(\calr^n,\mu))\\
&&\dim_{\calu(\calr,\mu)} H_k(\calr^n,\calu(\calr^n,\mu)) \otimes_{\calu(\calr^n,\mu)}
\calu(\calr,\mu)
\end{eqnarray*}

\subsection{Proof of Theorem~\ref{thm:betti-less normal subrelations}
and  ~\ref{thm:vanishing in degree plus one}}
 
Consider the strong extension 
\[1\rightarrow (\cals,\mu)\rightarrow (\calg,\mu)
\xrightarrow{p}(\calq,\nu) \to 1.\] 
The $E_1$-term of the corresponding spectral
sequence for is 
\begin{equation*}
E_1^{p,q}=\hom_{\calr(\calq,\nu)}\bigl( F_q(\calq), H^q(\cals;\calu(\calg,\mu))\bigr), 
\end{equation*}
where $F_\ast(\calq)$ is a projective 
resolution of $\calm(\calq^0)$ in $\compmod{\calr(\calq,\nu)}$. 

\smallskip

Thus, using Remark \ref{vanish}, $E_1^{p,q}=0$ for $0\le q\le d$. It follows 
that $E_r^{p,q}=0$ for $0\le q\le d$ and $1\le r\le\infty$. 
This implies that $E_r^{0,d+1}\cong E_\infty^{0,d+1}$
for $r\ge 2$ and $H^{d+1}(\calg,\calu(\calg,\mu)) \cong E_\infty^{0,d+1}$. 
Thus, 
\begin{equation}\label{eq:vanishing below degree d}
\betti_p(\calg)=0\text{ for $0\le p\le d$.}
\end{equation}
and 
\begin{align}\label{eq:formula for d-th betti}
\betti_{d+1}(\calg)
&=\dim_{\calu(\calg,\mu)}
  H^0\bigl(\calq;H^{d+1}(\cals;\calu(\calg,\mu)))\bigr)\\
&=\dim_{\calu(\calg,\mu)} \hom_{\calr(\calq,\nu)}
  \bigl(\calm(\calq^0), H^{d+1}(\cals;\calu(\calg,\mu))\bigr). \nonumber 
\end{align}
This completes the proof of
Theorem~\ref{thm:betti-less normal subrelations}.  For the proof of
Theorem \ref{thm:vanishing in degree plus one}, it remains to show
that $\betti_{d+1}(\calg)=0$ provided $\betti_{d+1}(\cals)<\infty$. 

\smallskip
Let us assume in addition that $(\calq,\nu)$ is an infinite equivalence relation.
Since $(\calq,\nu)$ is infinite, there
exists an infinite, amenable subrelation $\calr_{\hyp}\subset\calq$
(see~\citelist{\cite{gaboriau(2002b)}*{Proof of
    Th\'eor\`eme~6.8}\cite{cost}*{Proof of Proposition~III.3}}).  By
the same argument as for~(\ref{eq:formula for d-th betti}) applied to
$1\rightarrow\cals\rightarrow
p^{-1}(\calr_{\hyp})\rightarrow\calr_{\hyp}\rightarrow 1$ we get that
\begin{equation} \label{eq99} \betti_{d+1}\bigl(p^{-1}(\calr_{\hyp})\bigr)=
  \dim_{\calu(\calg,\mu)} \hom_{\calr(\calr_{\hyp},\nu)}
  \left(\calm(\calq^0) ,H^{d+1}(\cals;\calu(\calg,\mu))\right).
\end{equation}
 
Clearly, Equation \ref{eq99} implies
$\betti_{d+1}\bigl(p^{-1}(\calr_{\hyp})\bigr)\leq 
\betti_{d+1}(\cals)$, but this
holds also for the inverse image of any finite index sub-relation
$\calr'_{\hyp} \subset \calr_{\hyp}$. Hence, since such sub-relations
exist for any given index $n$, $$n \cdot
\betti_{d+1}(p^{-1}(\calr_{\hyp}))=
\betti_{d+1}(p^{-1}(\calr'_{\hyp}))\leq \betti_{d+1}(\cals)$$ for any
$n$, and so $\betti_{d+1}(p^{-1}(\calr_{\hyp})) =0$ as 
$\betti_{d+1}(\cals)<\infty$. 
Since the natural map
\[\hom_{\calr(\calq,\nu)}\left(\calm(\calq^0), H^{d+1}(\cals;\calu(\calg,\mu)\right) 
\rightarrow  
\hom_{\calr(\calr_{\hyp},\mu)}\left(
\calm(\calq^0),H^{d+1}(\cals;\calu(\calg,\mu)\right),\] 
is obviously injective, 
$\betti_{d+1}(\calg)=0$ follows now 
from~(\ref{eq:vanishing below degree d}) and~(\ref{eq:formula for d-th betti}). 

In the course of the proof, we assumed that $(\calq,\nu)$ was an equivalence relation. Again, the more general case of infinite discrete measured groupoids is a bit more 
tedious but analogous.

\subsection{Proof of Theorem~\ref{thm:vanishing by small ergodic dimension}}

Let us recall the definition of orbit equivalence. 

\begin{definition}
	We say that two countable 
	groups $\Gamma$ and $\Lambda$ are \emph{orbit equivalent} 
	if there is a probability space $(X,\mu)$ and essentially free $\mu$-preserving 
	actions of $\Gamma$ and $\Lambda$ such that for $\mu$-a.e.~$x\in X$ 
	we have $\Gamma x=\Lambda x$. Furthermore, $\Gamma$ and $\Lambda$ are 
	\emph{weakly orbit equivalent} if 
	there is a probability space $(X,\mu)$, essentially free $\mu$-preserving actions of 
	$\Gamma$ and $\Lambda$ on $X$, and Borel subsets $A\subset X$ 
	and $B\subset X$ such that $\Gamma A=X$ and $\Lambda B=X$ up to null sets and 
	such that for $\mu$-a.e.~$x\in X$ 
	we have $\Gamma x\cap A=\Lambda x\cap B$.
\end{definition}

\begin{lemma} \label{newlemma}
If the groups $Q_0$ and $Q_1$ are measure equivalent, then 
${\rm mcd}_{\bbC} (Q_0) \leq {\rm cd}_{\bbC}(Q_1)$. 
\end{lemma}
\begin{proof}
Furman~\cite{furman(1999a)} showed that groups are measure equivalent if and only 
if they are weakly orbit equivalent. 

Let us assume first that $Q_0$ and $Q_1$ are orbit equivalent. 
Let $(X,\mu)$ be a probability space equipped with m.p.~free actions of 
$Q_0$ and $Q_1$ such that the actions have the same orbit equivalence 
relation, which coincides with the translation groupoid. 
By~\cite{furman(1999a)}*{Lemma~2.2} 
one can assume that both actions are ergodic. 
So we have an identification of the corresponding groupoid rings 
$\calr(X\rtimes Q_0)\cong\calr(X\rtimes Q_1)$.  
Consider the following commutative diagram of functors: 
\[\xymatrix{  
\modules{\bbC Q_1}\ar[r] \ar@/_1cm/[rrdd]                     & 
\modules{L^\infty(X)\rtimes Q_1}\ar[r] & 
\compmod{L^\infty(X)\rtimes Q_1}\ar[d]^\cong\\ 
                                             &
                                             &
\compmod{\calr(X\rtimes Q_1)}\ar[d]^\cong        \\
                                                  &
                                                  &
\compmod{\calr(X\rtimes Q_0)}
}\]     
The horizontal arrows either denote completion or ring extension. The
vertical arrows are 
natural isomorphisms of abelian categories. Concerning the upper 
vertical arrow, we use that $L^\infty(X)\rtimes Q_1$ is dense 
in $\calr(X\rtimes Q_1)$, thus leading to an identification of 
complete $L^\infty(X)\rtimes Q_1$- with complete 
$\calr(X\rtimes Q_1)$-modules (see~\cite{thomrank}*{Lemma~4.4}). 
The first horizontal arrow is an exact functor 
because of $L^\infty(X)\rtimes Q_1\otimes_{\bbC Q_1}M\cong
L^\infty(X)\otimes_{\bbC}M$. The completion functors are exact 
and preserve projectives by~\cite{thomrank}*{Lemma~2.6 and
  Theorem~2.7}. 
Starting with a projective $\bbC Q_1$-resolution 
of $\bbC$ of length $n$, one obtains by following the arrows a 
projective resolution $P_\ast$ of $\calm(X)$ in $\compmod{\calr(X\rtimes
  Q_0)}$ of length $n$.

The general case where $Q_0$ and $Q_1$ are weak orbit equivalent 
demands only small modifications. For the upcoming discussion 
of full idempotents and Morita equivalence we 
recommend Lam's book~\cite{lam}*{Section~18}. Recall 
that an idempotent $p$ in a ring $R$ is called \emph{full} if 
the elements $rpr'$ with $r,r'\in R$ generate $R$ additively. 

\smallskip

Let $A\subset X$ and $B\subset X$ 
be Borel subsets of positive measure such that for a.e.~$x\in X$ we 
have $Q_0 x\cap A=Q_1 x\cap B$. In particular, this gives an
identification of the restricted translation groupoids 
\[\calr(X\rtimes Q_0)\vert_A\cong\calr(X\rtimes Q_1)\vert_B.\]
The characteristic functions 
$\chi_A$ and $\chi_B$ are \emph{full} idempotents in 
the groupoid rings $\calr(X\rtimes Q_0)$ and $\calr(X\rtimes Q_1)$, 
respectively. This is easily concluded from ergodicity and 
Lemma~\ref{lem:transitivity for ergodic groupoids}. 
Hence $\calr(X\rtimes Q_0)$ and 
\[\chi_A\calr(X\rtimes Q_0)\chi_A=\calr(X\rtimes Q_0\vert_A)\]
are Morita equivalent, i.e.~their module categories are 
equivalent as abelian categories. Explicitly, the Morita 
equivalence is given by $M\mapsto \chi_AM$. 
From this, it is obvious that the Morita equivalence 
restricts to the subcategories 
of complete modules. The same argument holds for $Q_1$ in place of $Q_0$. 
After replacing the lower vertical arrow in the diagram above 
by the equivalence of abelian categories
\[\xymatrix{
  \compmod{\calr(X\rtimes Q_1)}\ar[r]^{\simeq}\ar@{-->}[d] &
  \compmod{\calr(X\rtimes Q_1)\vert_B}\ar[d]^\simeq\\
  \compmod{\calr(X\rtimes Q_0)} &
  \compmod{\calr(X\rtimes Q_0)\vert_A}\ar[l]^\simeq
}\]
we can run the same kind of argument as before. 
This finishes the proof of the lemma.
\end{proof} 

We now turn to the prove of Theorem \ref{thm:vanishing by small ergodic dimension} in
the form needed for the application towards Theorem \ref{cor:app to hopf-singer}.
  
\begin{theorem} \label{secondvariant}
Let $1\rightarrow\Lambda\rightarrow\Gamma\rightarrow Q\rightarrow 1$
be a short exact sequence of groups. 
Suppose that $\betti_p(\Lambda)=0$ for $p>m$ and let $n={\rm mcd}_\bbC(Q)$. Then $\betti_p(\Gamma)=0$ for $p>n+m$. 
\end{theorem}  

\begin{proof}
Let $Q\curvearrowright (X,\nu)$ realize ${\rm mcd}_\bbC(Q)$. 
Consider Gaboriau's extension (see 
Subsection~\ref{subsec:example with groups by gaboriau}) 
of the form 
\begin{equation*}
1\rightarrow (Z\rtimes\Lambda,\mu\times\nu)\longrightarrow 
(Z\rtimes\Gamma,\mu\times\nu)\longrightarrow (X\rtimes Q,\nu)\rightarrow 1. 
\end{equation*}
The $E^1$-term of the corresponding spectral sequence is 
\begin{equation}\label{eq:E1-term}
E_1^{p,q}=\hom_{\calr(X\rtimes Q)}\bigl(P_p,
H^q(\calr(Z\rtimes\Lambda);\calu(Z\rtimes\Gamma))\bigr),   
\end{equation}
where $P_\ast$ is a projective resolution of $\calm(X)$ of length $n$ in the 
category $\compmod{\calr(X\rtimes Q)}$ (see Remark~\ref{rem:E1 term}). 
Hence $E_1^{p,q}=0$ whenever $p>n$. Moreover, $E_1^{p,q}=0$ if $q>m$, since
$H^q(\calr(Z\rtimes\Lambda);\calu(Z\rtimes\Gamma))=0$
for $q>m$, by Remark~\ref{vanish} and 
Equation~(\ref{eq:betti identification}) . This yields the theorem 
since $E_{\infty}^{p,q}=0$ for all
$p+q > n+m$, and thus $\betti_i(\Gamma) =\betti_i(Z\rtimes\Gamma) =0$ for $i>n+m$. 
\end{proof}

Clearly, Theorem \ref{secondvariant} and Lemma \ref{newlemma} imply Theorem \ref{thm:vanishing by small ergodic dimension}.

\end{document}